\documentclass[12pt]{iopart}

 \expandafter\let\csname equation*\endcsname\relax 
 \expandafter\let\csname endequation*\endcsname\relax 
\usepackage{amsmath}
\usepackage{amsthm}		
\usepackage{amssymb} 
\usepackage{graphicx}
\usepackage{xcolor}
\usepackage[colorlinks=true, linkcolor=blue]{hyperref}
\usepackage{cite}

\hypersetup{
pdftitle={Boundary value problems and singularity theory},
pdfauthor={Christian Offen},
pdfsubject={Bifurcation Theory, symplectic geometry, boundary value },
pdfkeywords={bifurcation, boundary value problems, symplectic geometry, Lagrangian, Hamiltonian systems, integrable systems, periodic orbits, invariant tori, pitchfork, saddle-node, singularity theory},
bookmarksdepth=10,
}

\def\d{\mathrm{d}}
\def\p{\partial }

\def\D{\mathrm{D}}

\def\Id{\mathrm{Id}}
\def\e{\epsilon}

\def\Z{\mathbb{Z}}

\def\R{\mathbb{R}}

\newtheorem{definition}{Definition}[section]
\newtheorem{theorem}{Theorem}[section]
\newtheorem{prop}{Proposition}[section]
\newtheorem{corollary}{Corollary}[section]
\newtheorem{lemma}[theorem]{Lemma}
\newtheorem{remark}{Remark}[section]

\begin{document}

\title{Bifurcation of solutions to Hamiltonian boundary value problems}
\author{R I McLachlan$^1$ and C Offen$^1$}
\address{$^1$ Institute of Fundamental Sciences, Massey University, Palmerston North, New Zealand}
\eads{\mailto{r.mclachlan@massey.ac.nz}, \mailto{c.offen@massey.ac.nz}}

\begin{abstract}
A bifurcation is a qualitative change in a family of solutions to an equation produced by varying parameters. In contrast to the local bifurcations of dynamical systems that are often related to a change in the number or stability of equilibria, bifurcations of boundary value problems are global in nature and may not be related to any obvious change in dynamical behaviour. Catastrophe theory is a well-developed framework which studies the bifurcations of critical points of functions. In this paper we study the bifurcations of solutions of boundary-value problems for symplectic maps, using the language of (finite-dimensional) singularity theory. We associate certain such problems with a geometric picture involving the intersection of Lagrangian submanifolds, and hence with the critical points of a suitable generating function. Within this framework, we then study the effect of three special cases: (i) some common boundary conditions, such as Dirichlet boundary conditions for second-order systems, restrict the possible types of bifurcations (for example, in generic planar systems only the A-series beginning with folds and cusps can occur); (ii) integrable systems, such as planar Hamiltonian systems, can exhibit a novel {\em periodic pitchfork} bifurcation; and (iii) systems with Hamiltonian symmetries or reversing symmetries can exhibit restricted bifurcations associated with the symmetry. This approach offers an alternative to the analysis of critical points in function spaces, typically used in the study of bifurcation of variational problems, and opens the way to the detection of more exotic bifurcations than the simple folds and cusps that are often found in examples. 
\end{abstract}



\submitto{\NL}

\maketitle

\section{Introduction}\label{subsec:motivation}


Let $x,y$ be the standard coordinates of $\R^2$ with the standard symplectic structure $\omega = \d x \wedge \d y$.
We consider the two-parameter family of Hamiltonian systems defined by
\[
H_\mu \colon \R^2 \to \R, \quad 
H_\mu(x,y) = y^2+\mu_1 x + \mu_2 x^2+x^4.
\]
For $\mu \in \R^2$ the time-$\tau$-flow map $\phi_\mu \colon \R^2 \to \R^2$ of the Hamiltonian vector field $X_H$ defined by $\d H = \omega(X_H,.)$ assigns to initial values $(x(0),y(0))$ the solution of Hamilton's equation
\begin{equation}\label{eq:HamSys}
\begin{cases}
\dot x (t) &=  \frac{\p H_\mu}{\p y} (x(t),y(t))\\
\dot y (t) &= -\frac{\p H_\mu}{\p x} (x(t),y(t))
\end{cases}
\end{equation}
at time $\tau$. Each map $\phi_\mu$ is symplectic, i.e.\ $\phi_\mu^\ast \omega = \omega$.
Let us consider the Dirichlet boundary value problem
\begin{equation}\label{eq:bd-problem-example}
x(0)=x^\ast
\qquad
x(\tau)=X^\ast
\end{equation}
for $x^\ast,X^\ast \in \R$.
In other words, we look for orbits of the Hamiltonian flow which start on the line $\{x^\ast\}\times \R$ in the phase space $\R^2$ and end on the line $\{X^\ast\}\times \R$ at time $\tau$. 
Since solutions to initial value problems are unique, we can specify a solution to \eqref{eq:HamSys} and \eqref{eq:bd-problem-example} by the value $y(0)$. This means we seek solutions $y(0)$ to the equation
\begin{equation}\label{eq:bd-problem-example-exact}
(x\circ\phi_\mu)(x^\ast,y(0)) = X^\ast.
\end{equation}

\begin{figure}
\begin{center}
\includegraphics[width=0.8\textwidth]{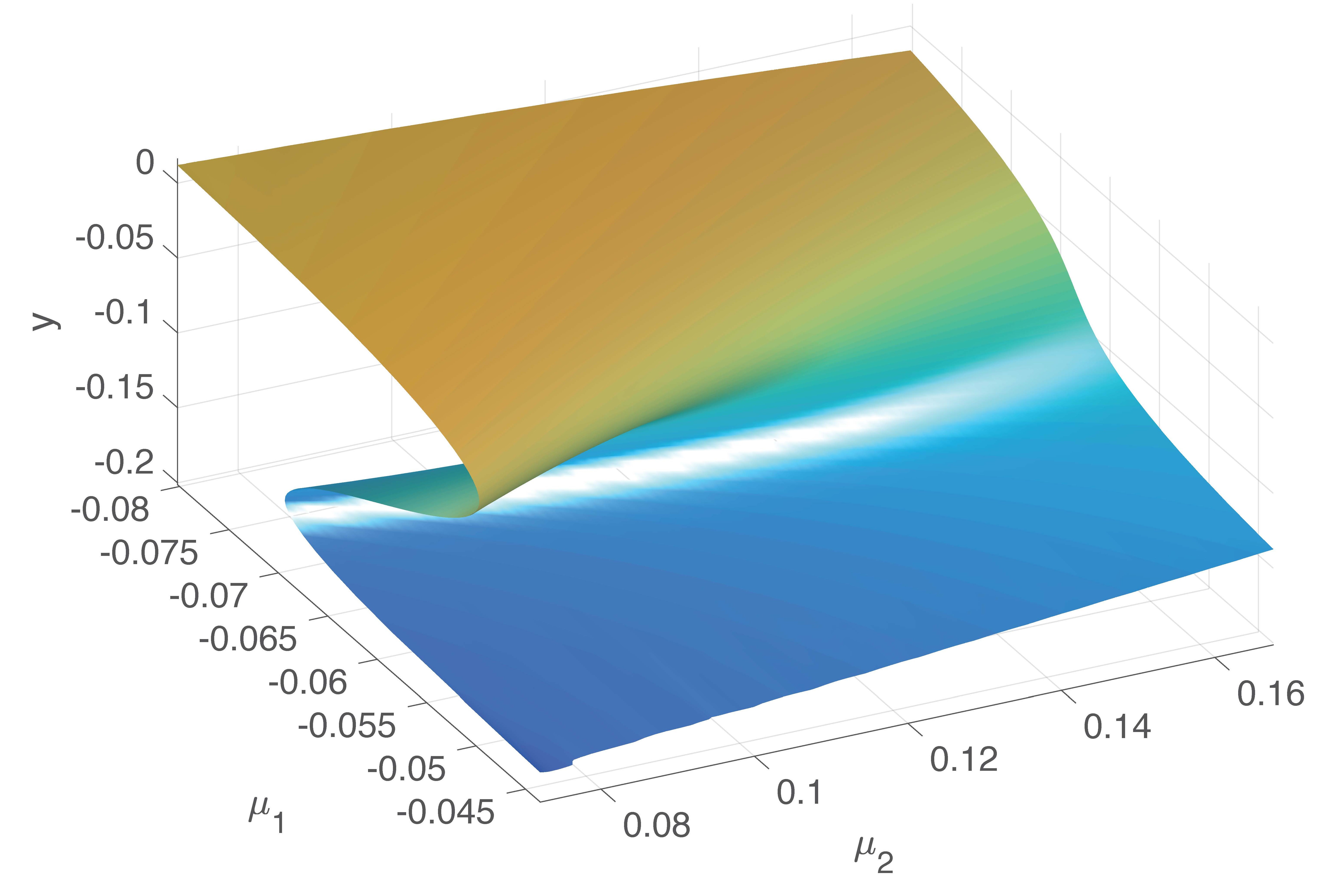}
\end{center}
\caption{A plot of solutions $y$ to \eqref{eq:bd-problem-example-exact} over the parameter space for $\tau=4$, $x^\ast = 0.2 = X^\ast$ shows a cusp bifurcation. The Hamiltonian flow is approximated using the symplectic leapfrog method with step size $0.1$.}\label{fig:cusp-example-problem}
\end{figure}

Figure \ref{fig:cusp-example-problem} shows a bifurcation diagram, i.e.\ a plot of solutions $y(0)$ of \eqref{eq:bd-problem-example-exact} over the $\mu$-plane for the integration time $\tau=4$ and the boundary values $x^\ast = 0.2 = X^\ast$.
The plot shows a \textit{cusp bifurcation}. The bifurcation is classically known as one of the seven elementary catastrophes classified by Thom \cite{ThomOriginal}. 
(See e.g. \cite[\S 6]{Gilmore1993catastrophe} for a detailed, elementary description of the elementary catastrophes).
In the context of catastrophe theory the bifurcation diagram is referred to as a catastrophe set.
It persists under small perturbations, i.e.\ a small perturbation of the Hamiltonian results in a bifurcation diagram which qualitatively looks like\footnote{More precisely, the perturbed and unperturbed bifurcation diagrams correspond to unfoldings of smooth, scalar valued maps which are equivalent in the left-right category of unfoldings \cite[p.68ff]{lu1976singularity}.} the bifurcation diagram of the unperturbed system shown in figure \ref{fig:cusp-example-problem}.
This means that the above bifurcation is a \textit{generic} phenomenon in Dirichlet problems for two parameter families  of Hamiltonian systems. \\

The occurrence of a cusp bifurcation is related to the symplecticity of the maps $\phi_\mu$. Indeed, it is an instance of the fact that many boundary value problems for symplectic maps are governed by catastrophe theory which we will explain rigorously in sections \ref{subsec:connection-SingularityTheory} and \ref{sec:bd-problemssymplmaps}.
Let us sketch the connection here for Dirichlet problems for symplectic maps on $\R^{2n}$. Let $\phi_\mu \colon (\R^{2n},\omega) \to (\R^{2n},\omega)$ denote a family of symplectic maps.
We equip $(\R^{2n},\omega)$ with symplectic coordinates $x^1,\ldots,x^n,y_1,\ldots,y_n$ such that $\omega= \sum_{j=1}^n \d x^j \wedge \d y_j$.
Let $\pi_1,\pi_2 \colon \R^{2n}\times \R^{2n} \to \R^{2n}$ denote the projection to the first or second factor of the cross product $\R^{2n}\times \R^{2n}$, respectively.
We equip $\R^{2n}\times \R^{2n}$ with the symplectic form $\omega \oplus (-\omega) := \pi_1^\ast \omega - \pi_2^\ast \omega$ and obtain coordinates $x^1,\ldots,x^n,y_1,\ldots y_n, X^1,\ldots,X^n,Y_1,\ldots Y_n$ on $\R^{2n}\times \R^{2n}$ by pulling back the coordinates $x^1,\ldots,x^n,y_1,\ldots,y_n$ from $\R^{2n}$ with $\pi_1$ and $\pi_2$.
The graphs of the symplectic maps can be embedded as Lagrangian submanifolds:
\begin{equation*}
\Gamma_\mu = \{ (x,y,\phi_\mu(x,y)) \, | \, (x,y)\in \R^{2n}\} \subset (\R^{2n}\times \R^{2n}, \omega \oplus (-\omega)).
\end{equation*}
The 1-form
\[
\alpha 
= \sum_{j=1}^n \left((x_j - x^\ast )\d y_j - (X_j - X^\ast )\d Y_j\right)
\]
on $\R^{2n}\times \R^{2n}$ is a primitive of $\omega \oplus (-\omega)$ and, therefore, closed on the simply connected submanifolds $\Gamma_\mu$ for each $\mu$.
Thus, there exists a family of primitives $S_\mu \colon \Gamma_\mu \to \R$ with $\d S_\mu = \iota_\mu^\ast\alpha$, where $\iota_\mu \colon \Gamma_\mu \hookrightarrow \R^{2n}\times\R^{2n}$ is the natural inclusion.
Where $y_j,Y_j$ ($1\le j \le n$) constitutes a local coordinate system on $\Gamma_\mu$, i.e.\ where $\det \left( \frac{\p (y_j \circ \phi_\mu)}{\p x^i}\right)_{i,j} \not =0$, the problem
\[
(x \circ \phi_\mu)(x^\ast,y)=X^\ast
\]
is equivalent to $\iota_\mu^\ast \alpha=0$ or
\[
\d S_\mu =0.
\]

%
%
%
We can conclude that the bifurcations in the boundary value problem \eqref{eq:bd-problem-example} for families of symplectic maps $\phi_\mu$ behave like the \textit{gradient-zero-problem}, which is well understood and known under the name \textit{catastrophe theory}.\\

Bifurcations of the problem $\d S_\mu =0$ were first classified by Thom \cite{ThomOriginal}. All singularities that occur generically in gradient-zero-problems with at most 4 parameters are stably right-left equivalent to one of the following seven elementary catastrophes:

\begin{center}
  \begin{tabular}{| c | l | c | r | c |}
    \hline
	ADE class & name & germ & truncated miniversal deformation  \\
    \hline
    \hline
$A_2$ &fold & $x^3$ & $x^3+ \mu_1 x$ \\
$A_3$ &cusp & $x^4$ & $x^4+\mu_2 x^2+ \mu_1 x$ \\
$A_4$ &swallowtail & $x^5$ & $x^5+\mu_3 x^3+\mu_2 x^2+ \mu_1 x$  \\
$A_5$ &butterfly & $x^6$ & $x^6+\mu_4 x^4+\mu_3 x^3+\mu_2 x^2+ \mu_1 x$  \\
\hline 
$D_4^+$ &hyperbolic umbilic & $x^3+xy^2$ & $ x^3+xy^2 + \mu_3 (x^2-y^2) +\mu_2 y + \mu_1 x$  \\
$D_4^-$ &elliptic umbilic & $x^3-xy^2$ & $ x^3-xy^2 + \mu_3 (x^2+y^2) +\mu_2 y + \mu_1 x$  \\
$D_5$ &parabolic umbilic & $x^2y+y^4$ & $ x^2y+y^4 +\mu_4 x^2  + \mu_3 y^2 +\mu_2 y + \mu_1 x$.  \\
 \hline
  \end{tabular}\label{tabel:elementaryCats}
\end{center}

Thom's results were extended by Arnold in a series of papers including \cite{ArnoldNormal1,ArnoldNormalUni,ArnoldNormal2,ArnoldNormal3}. Results are summarised in \cite{Arnold1,Arnold2}. 
A good introduction to the algebraic framework of catastrophe theory and the results obtained by Thom is \cite{lu1976singularity}. For a detailed, elementary consideration of singularities with low multiplicity we refer to \cite{Gilmore1993catastrophe} where the reader will also find examples of applications in physics.\\

In the setting used in the example, Dirichlet boundary conditions $x=x^\ast$, $X=X^\ast$ can be represented as Lagrangian submanifolds
\[
\Lambda = \{ (x,y,X,Y)\, |\, x=x^\ast, X=X^\ast, y,Y \in \R^{n} \}
\]
of $(\R^{2n}\times \R^{2n},\omega \oplus (-\omega))$. Solutions of the boundary value problem then correspond to elements of the Lagrangian intersection $\Lambda \cap \Gamma_\mu$.
We can represent other boundary conditions as well: Periodic boundary conditions $(X,Y)=(x,y)$ 
are described by the submanifold
\[
\Lambda = \{ (x,y,X,Y)\, |\, X=x,Y=y, x,y \in \R^{n} \}
\]
which is also Lagrangian. Neumann boundary conditions are conditions on the derivatives of second order ODEs at the boundaries of a time interval. In a first order formulation this often translates to conditions on the $y$-variables, $y=y^\ast$, $Y=Y^\ast$ analogous to Dirichlet problems. The representing submanifold
\[
\Lambda = \{ (x,y,X,Y)\, |\,y=y^\ast,Y=Y^\ast, x,X \in \R^{n} \}
\]
is Lagrangian. More generally, linear boundary conditions
\[
A \begin{pmatrix}
x\\ y\\ X \\Y
\end{pmatrix}
=\begin{pmatrix}
\alpha\\ \beta
\end{pmatrix}
\]
with $A \in \R^{2n \times 4n}$ of maximal rank and $\alpha, \beta \in \R^{n}$ constitute Lagrangian submanifolds if and only if $V^t J V = 0$ where the columns of $V\in \R^{4n \times 2n}$ span the kernel of $A$ and $J$ represents the symplectic form $\omega \oplus (-\omega)$ on $\R^{2n}\times \R^{2n}$, i.e.\
\[
J=\begin{pmatrix}
0&-I_n&0&0\\
I_n&0&0&0\\
0&0&0&I_n\\
0&0&-I_n&0
\end{pmatrix}.
\]
Here $I_n \in \R^{n \times n}$ denotes the identity matrix. If $n=1$ then
\[
x+y=\alpha, \quad X+Y=\beta
\quad \text{corresponding to}\quad
A = \begin{pmatrix}
1&1&0&0\\0&0&1&1
\end{pmatrix},\,
V=\begin{pmatrix}
1&1\\-1&-1\\1&-1\\-1&1
\end{pmatrix}
\]
is an example of a Lagrangian boundary condition while
\[
x-y+X-Y = \alpha,\quad
x+y+X+Y=\beta
\quad \text{corresponding to}\]
\[
A = \begin{pmatrix}
1&-1&1&-1\\1&1&1&1
\end{pmatrix},\,
V=\begin{pmatrix}
1&0\\0&1\\-1&0\\0&-1
\end{pmatrix}
\]
is an example of a boundary condition that is not Lagrangian. We will show that the bifurcation behaviour of Lagrangian boundary value problems are governed by catastrophe theory (gradient-zero-problem) while the bifurcation of non-Lagrangian boundary value problems in general position behave like the root-of-a-function problem / equilibria-of-vector-fields-problem: we can expect the same bifurcation behaviour as in the problem
\[
F_\mu(x,y)=0
\]
for a family of maps $F_\mu \colon \R^{2n} \to \R^{2n}$. Here, in contrast to Lagrangian boundary value problems, the maps $F_\mu$ are not necessarily gradients of scalar valued maps. See \cite{TopoStability} for a treatment of singularities occurring in the zeros-of-a-function problem $F_\mu(z)=0$, $F_\mu \colon \R^m \to \R^k$. An alternative, earlier reference are the lectures \cite{walllectures} which are based on \cite{MatherI,MatherII,MatherIII,MatherIV,MatherV,MatherVI}.\\

The general idea to translate boundary value problems of symplectic maps to a catastrophe theory setting will be explained in the following two sections and can be summarised as follows:

\begin{itemize}
\item
Symplectic maps and Lagrangian boundary conditions constitute Lagrangian submanifolds in a cotangent bundle.
\item
We will show in lemma \ref{lem:graphoverzero} that
after shrinking the involved manifolds and the parameter space and after applying a symplectic transformation all involved Lagrangian submanifolds can assumed to be graphs over the zero-section.
\item
The Lagrangian submanifolds can be written as images of exact 1-forms and the primitives of the forms coming from the boundary conditions can be substracted from the primitives related to the symplectic maps.
\item
We obtain a family of smooth functions whose critical points correspond to solutions of the boundary value problem.
\end{itemize}

For the motivational example this means that Dirichlet boundary value problems (making use of global structure of the phase space) turn out to admit a translation to a \textit{local} problem involving intersections of Lagrangian submanifolds. This allows a treatment with basic tools from symplectic geometry making use of the Darboux theorem that locally all symplectic manifolds are symplectomorphic. In contrast, the treatment of global Lagrangian intersection problems crucially involves the topology of the manifolds. Let us refer at this point to Arnold's conjecture, which (in a special case which has been proved) asserts a lower bound for the amount of fixed points of Hamiltonian diffeomorphisms. The statement can be interpreted as a lower bound for the number of intersection points of a Lagrangian submanifold and a Hamiltonian-isotopic Lagrangian submanifold 
\cite{fukaya2010lagrangian}. Generalisations of Arnold's conjecture can be found in \cite{Oh1995,Ciriza2000,Polterovich2001}.\\

In \cite{WeinstBvPHam} Weinstein develops a similar geometric picture which he uses to address global existence questions about solutions to boundary value problems. The focus is on fixed boundary value problems in Hamiltonian systems with regular (i.e., nonsingular and nonbifurcating) manifolds as solutions. In contrast to \cite{WeinstBvPHam}, the solutions considered in this paper are typically isolated and we consider families of boundary value problems and local bifurcations of their solutions. In addition, in our work time is either fixed or interpreted as a parameter, while in \cite{WeinstBvPHam} a fixed energy-time relation is imposed. 
\\

The geometric picture of intersecting Lagrangian submanifolds in this paper and in \cite{WeinstBvPHam} is to be contrasted with studies that involve Lagrangian intersections in the phase space. In this paper a pair of Lagrangian submanifolds represents a boundary value problem. The manifolds are submanifolds of the product of the phase space with itself. In contrast, Lagrangian submanifolds of phase spaces often represent objects of dynamical interest (e.g.\ invariant manifolds) and their intersections provide information about the dynamics; see e.g.\ \cite{Haro,LMS,LMR}.
\\

There is another approach to the study of bifurcations of boundary value problems, which we mention briefly to contrast with that of the present paper. It is often used for PDEs. 
Indeed, many weak formulations of PDEs arise as variations of functionals such that critical points of functionals correspond to weak solutions of the PDE. Let us consider the following example which is a generalisation of the much-studied Bratu problem. (See \cite{Mohsen201426} for a review of the Bratu problem.)\\

\begin{figure}
\begin{center}
\includegraphics[width=0.6\textwidth]{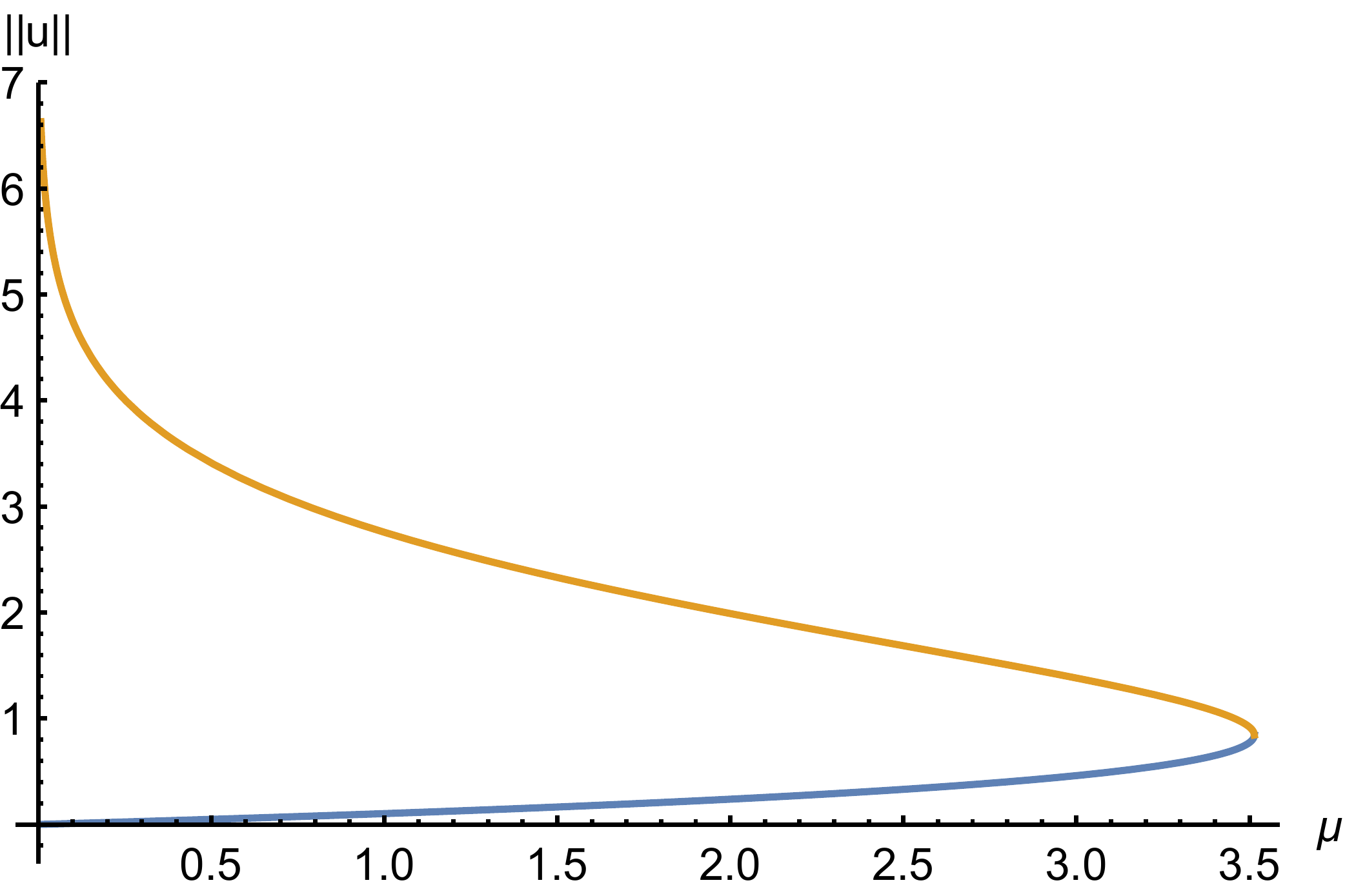}
\end{center}
\caption{The bifurcation diagram of the 1-parameter, 1-dimensional Bratu problem shows a fold bifurcation at $\mu\approx 3.5138$. The diagram was obtained using analytic solutions.\cite{Mohsen201426} }\label{fig:bratu1bifur}
\end{figure}

For a family of smooth functions $f_\mu \colon \R \to \R$ and the cube $\Omega = (0,1)^d \subset \R^d$ we consider the following PDE with Dirichlet boundary conditions
\begin{equation}\label{eq:PoissonPDE}
\begin{cases}
\Delta u &= f_\mu'(u)\\
u|_{\partial \Omega} &= 0.
\end{cases}
\end{equation}
In the classical 1-parameter, 1 dimensional Bratu problem with $f_\mu(u)=-\mu e^u$ the solutions undergo a bifurcation from infinity at $\mu=0$ and a fold bifurcation at $\mu\approx 3.5138$, see figure \ref{fig:bratu1bifur}.
The weak formulation of \eqref{eq:PoissonPDE} is given as
\begin{equation}\label{eq:weekform}
\forall v \in H : 
\quad \int_{\Omega} \left( \langle \nabla u (x), \nabla v (x) \rangle + f_\mu'(u(x)) v(x) \right)\d x=0,
\end{equation}
where $x=(x_1,\ldots,x_d)$ are coordinates on $\R^d$, $\langle.,. \rangle$ denotes the euclidean inner product on $\R^d$ and $H=H_0^1$ is the closure in the Sobolev space $W^{1,2}(\Omega)$ of the space $\mathcal C^\infty_c(\Omega)$ of smooth functions with compact support on the open cube. Notice that $H$ is a Hilbert space with inner product induced from $W^{1,2}(\Omega)$.
The equation \eqref{eq:weekform} can be written as
\[\forall v \in H :  D S|_u(v)=0,\]
where $D$ denotes the Fr\'echet derivative and $S\in H^\ast$ is the functional
\[ S_\mu(u) = \int_\Omega \left( \frac 12 \langle \nabla u, \nabla u \rangle + f_\mu(u)\right) \d x. \]

Bifurcation points $(\mu^\ast,u^\ast)$ of solutions to \eqref{eq:weekform} thus correspond to zeros of the Fr\'echet derivative of $S$.\\

There exist general statements about basic bifurcations in the critical points of functional problem, typically under technical assumptions which allow Lyapunov-Schmidt reduction to a finite dimensional problem, see e.g.\cite{Kielhoefer2012} or \cite{MR3161608}. Note that even in the case $d=1$, i.e.\ boundary value problems for ODEs, this approach uses the setting of functional analysis whereas ours is purely finite-dimensional throughout.\\

The paper is organised as follows: in section \ref{subsec:connection-SingularityTheory} we introduce a setting to analyse local intersection problems for two families of Lagrangian submanifolds.
Although classically known as the morsification of Lagrangian intersections, we present the procedure of translating local Lagrangian intersection problems to the problem of finding critical points of a scalar valued map (gradient-zero-problem) in detail because an accurate picture is needed for the statements in section \ref{sec:bd-problemssymplmaps}.\\


In section \ref{sec:bd-problemssymplmaps} we use this geometric framework to obtain the following results:

\begin{itemize}
\item
Bifurcations in Lagrangian boundary value problems for families of symplectic maps are governed by catastrophe theory.
\item
Some boundary conditions restrict the class of bifurcations that can occur. We will show how to test boundary conditions on whether and how they impose restrictions.
For example, in planar Dirichlet problems only $A$-series bifurcations can occur.
\item
Extra structure in the phase space can lead to new bifurcations not captured by catastrophe theory. We develop a detailed model showing that a pitchfork bifurcation occurs generically in 1-parameter families of completely integrable systems, for example in planar Hamiltonian systems. As it is related to the generic occurrence of periodic orbits, we refer to this bifurcation as \textit{periodic pitchfork bifurcation}. 
\item
We explain how symmetries and reversing symmetries of symplectic maps relate to symmetries of the corresponding gradient-zero-problem.
\end{itemize}


\section{Lagrangian intersections and catastrophe theory}\label{subsec:connection-SingularityTheory}

In the 
introduction we have seen that bifurcations of Dirichlet Hamiltonian boundary value problems are governed by catastrophe theory. Indeed, this is an instance of a more general phenomenon.
We will present the details of a theory in which certain types of Hamiltonian boundary value problems are systematically translated into a geometric picture involving the intersection of Lagrangian submanifolds, and then into a gradient zero problem. The process is classically known as \textit{morsification of Lagrangian intersections} \cite[p.33]{LagrangianInterTheo}. 
Bifurcations in boundary value problems will be shown to correspond to singularities of families of scalar valued functions.
Introductory references to symplectic geometry are \cite{Libermann1987, SymplecticTopology}. A classical reference for properties of Lagrangian submanifolds is \cite{WeinstLagr}.\\




Let $N$ be a smooth manifold. We equip the cotangent bundle $\pi \colon T^\ast N \to N$ with the symplectic structure $\omega = \d \lambda$, where $\lambda$ is the Liouvillian-1-form on $T^\ast N$ which is canonically defined by
\begin{equation}\label{eq:Liouvillian}
\lambda_\alpha(v) = \alpha(\d \pi|_{\alpha}(v)) \quad \text{for all }\, \alpha \in T^\ast N, v \in T_\alpha T^\ast N.
\end{equation}

\begin{definition}[local Lagrangian intersection problem]
Let $\mu^\ast \in \R^p$ and let $I$ be an open neighbourhood of $\mu^\ast$. Consider two smooth families $(\Lambda_\mu)_{\mu \in  I}$ and $(\Gamma_\mu)_{\mu \in  I}$ of Lagrangian submanifolds of $(T^\ast N,\omega)$. Let $z^\ast \in \Lambda_{\mu^\ast} \cap \Gamma_{\mu^\ast}$. The collection $\big( (\Lambda_\mu)_{\mu \in  I}, (\Gamma_\mu)_{\mu \in  I},(\mu^\ast,z^\ast) \big) $ is called a local Lagrangian intersection problem in $(T^\ast N,\omega)$.
\end{definition}

The following lemma asserts that two families of Lagrangian submanifolds can locally be mapped symplectomorphically to graphical Lagrangian families, i.e.\ without loss of generality all submanifolds can assumed to be graphs over the zero-section. 

\begin{lemma}\label{lem:graphoverzero}
Let $\big( (\Lambda_\mu)_{\mu \in \tilde I}, (\Gamma_\mu)_{\mu \in  \tilde I},(\mu^\ast,z^\ast) \big) $ be a local Lagrangian intersection problem in $(T^\ast N,\omega)$.
There exists 
an open neighbourhood $V\subset T^\ast N$ of $z^\ast$, 
a symplectomorphism $\Psi$ defined on $T^\ast (\pi(V))$ fixing $z^\ast$ and an open neighbourhood $I\subset \tilde I \subset \R^p$ of $\mu^\ast$ such that for all $\mu \in I$ the Lagrangian submanifolds $\Psi(\Lambda_\mu\cap V)$ and $\Psi(\Gamma_\mu\cap V)$ are graphical in $T^\ast (\pi(V))$, i.e.\ $\pi|_{\Psi(\Gamma_\mu \cap V)}$ and $\pi|_{\Psi(\Lambda_\mu \cap V)}$ are injective.
\end{lemma}
 

\textit{Proof.} Let us shrink $N$ (and $\Lambda_\mu,\Gamma_\mu\subset T^\ast N$ accordingly) to a coordinate neighbourhood of $\pi(z^\ast)\in N$ with coordinates $x^1,\ldots,x^n$.
Define coordinates $q^1,\ldots,q^n,p_1,\ldots,p_n$ on $T^\ast N$ centred at $z^\ast$ by

\begin{align*}
q^i \colon T^\ast N \to \R \quad &q^i = x^i \circ \pi \\
p_i\colon T^\ast N \to \R \quad &\gamma \mapsto \gamma\left(\left.\frac{\p}{\p x^i}\right|_{\pi(\gamma)}\right)-z^\ast\left(\left.\frac{\p}{\p x^i}\right|_{\pi(z^\ast)}\right)
\end{align*} for $i \in \{1,\ldots,n\}$.\footnote{If $T^\ast N=T^\ast\R^n$ with standard coordinates $(x^i,y_i)$ and $z^*=(x^{*i},y^*_i)$ then $q^i=x^i$, $p_i=y_i - y_i^*$.}
We have $\lambda = \sum_{i=1}^n p_i \d q^i$ and $\omega = \sum_{i=1}^n \d p_i \wedge \d q^i$.
Let $n_1 = \dim( T_{z^\ast}\Lambda_{\mu^\ast} \cap T_{z^\ast}\Gamma_{\mu^\ast})$ and $n_2=n-n_1$. Consider the symplectic basis $\left.\frac{\p}{\p q^1}\right|_{z^\ast},\ldots, \left.\frac{\p}{\p q^n}\right|_{z^\ast},\left.\frac{\p}{\p p_1}\right|_{z^\ast},\ldots,\left.\frac{\p}{\p p_n}\right|_{z^\ast}$ on $T_{z^\ast}T^\ast N$. In \cite{coisoPair} all pairs of coisotropic linear subspaces of finite dimensional vector spaces are classified. Using this classification result, there exists a linear symplectic map on $T_{z^\ast}T^\ast N$ represented by a matrix $M_0$ mapping $T_{z^\ast}\Lambda_{\mu^\ast}$ to the Lagrangian subspace spanned by the columns of $A_0$ and $T_{z^\ast}\Gamma_{\mu^\ast}$ to the space spanned by the columns of $B_0$ where
\[
A_0=\begin{pmatrix}
\Id_{n_1}&0_{n_1\times n_2}\\
0_{n_2\times n_1}&\Id_{n_2}\\
0_{n_1\times n_1}&0_{n_1\times n_2}\\
0_{n_2\times n_1}&0_{n_2\times n_2}
\end{pmatrix}
\quad \text{and} \quad
B_0=\begin{pmatrix}
\Id_{n_1}&0_{n_1\times n_2}\\
0_{n_2\times n_1}&0_{n_2\times n_2}\\
0_{n_1\times n_1}&0_{n_1\times n_2}\\
0_{n_2\times n_1}&\Id_{n_2}
\end{pmatrix}.
\]
Here $\Id_{k}$ denotes the $k$-dimensional identity matrix and $0_{k \times l}$ the zero matrix in $\R^{k \times l}$. Using the linear, symplectic transformation represented by the matrix
\[
M =
\begin{pmatrix}
\Id_{n}&\Id_{n}\\
0_{n\times n}&\Id_{n}
\end{pmatrix}
\]
$A_0$ and $B_0$ are mapped to matrices of the form $\begin{pmatrix}\Id_{n}\\ \ast\end{pmatrix}$. Their columns span Lagrangian subspaces which are graphs over $\mathrm{span}\left\{ \left.\frac{\p}{\p q^j}\right|_{z^\ast} \right\}_{1\le j \le n}$. Using the coordinate system $q^1,\ldots,q^n,p_1,\ldots,p_n$, the matrix $M \cdot M_0$ defines a symplectic map $\Psi$ on $T^\ast N$.
By construction, $\Psi(z^\ast)=z^\ast$ and for all $\mu$ near $\mu^\ast$ the Lagrangian manifolds $\Psi(\Lambda_\mu)$ and $\Psi(\Gamma_\mu)$ are locally around $z^\ast$ graphs over the zero-section in $T^\ast N$.\qed

\begin{remark}\label{rem:zerosecinv}
If we first apply a parameter-dependent, symplectic change of coordinates to the local Lagrangian intersection problem $\big( (\Lambda_\mu)_{\mu \in \tilde I}, (\Gamma_\mu)_{\mu \in  \tilde I},(\mu^\ast,z^\ast) \big) $ such that $\Lambda_\mu$ lies in the zero-section then the proof of lemma \ref{lem:graphoverzero} provides a symplectic map leaving the zero section invariant. Two symplectomorphisms mapping $\Lambda_\mu$ to the zero-section and $z^\ast$ to the same point $x$ differ only by a symplectomorphism leaving the zero-section invariant (locally around $x$).
\end{remark}

\begin{remark}\label{rem:tangentgraphical}
The lemma shows that we can transform the intersection problem such that not only the families of Lagrangian submanifolds become graphical but also their tangent spaces.
\end{remark}

Let us recall the notion of a catastrophe set in catastrophe theory.

\begin{definition}[catastrophe set of a singularity]
Let $N$ be a smooth manifold and let $(h_\mu)_{\mu\in I}$ be a family of smooth maps $h_\mu \colon N \to \R$. The set
\[\{ (\mu, x)\in I \times N \; | \; \nabla h_\mu(x)=0 \}\]
is called the \textit{catastophe set of the family $(h_\mu)_{\mu\in I}$}.
\end{definition}

We may use a similar definition in the setting of local Lagrangian intersection problems:

\begin{definition}[catastrophe set of intersection problems]
Let $\mathfrak L = \big( (\Lambda_\mu)_{\mu \in I}, (\Gamma_\mu)_{\mu \in  I},(\mu^\ast,z^\ast) \big) $ be a local Lagrangian intersection problem in $(T^\ast N,\omega)$ such that $\pi|_{\Lambda_\mu}$ and $\pi|_{\Gamma_\mu}$ are injective for all $\mu \in I$. The set
\[\{ (\mu,\pi(z))\; |\; \mu \in I,\, z \in \Lambda_{\mu} \cap \Gamma_{\mu} \}\] is called catastrophe set of the local Lagrangian intersection problem $\mathfrak L$.
\end{definition}

The following theorem asserts that the intersection of two families of Lagrangian submanifolds locally behaves like the gradient-zero problem for a smooth family of local maps on $N$.

\begin{theorem}[Singularity-Lagrangian linkage theorem]\label{thm:linkageThm}
Let $\mathfrak L = \big( (\Lambda_\mu)_{\mu \in \tilde I}, (\Gamma_\mu)_{\mu \in \tilde I},(\mu^\ast,z^\ast) \big) $ be a local Lagrangian intersection problem in $(T^\ast N,\omega)$.
There exists an open neighbourhood $V \subset T^\ast N$ of $z^\ast$, a smooth family $(h_\mu)_{\mu \in I}$ of smooth maps $h_\mu \colon \pi(V) \to \R$ 
and a symplectomorphism $\Psi$ on $T^\ast (\pi(V))$ fixing $z^\ast$ and mapping $\mathfrak L$ to a local Lagrangian intersection problem
\[
\big( (\Psi(\Lambda_\mu \cap V) )_{\mu \in I}, (\Psi(\Gamma_\mu\cap V))_{\mu \in I},(\mu^\ast,z^\ast) \big)
\]
in $T^\ast \pi(V)$ with the same catastrophe set as $(h_\mu)_{\mu\in I}$.\\

In other words, up to a local symplectomorphism for each local Lagrangian intersection problem there exists a smooth family of smooth maps with the same catastrophe set.

\end{theorem}


\textit{Proof.} By lemma \ref{lem:graphoverzero} there exists 
an open neighbourhood $V\subset T^\ast N$ of $z^\ast$, 
a symplectomorphism $\Psi$ defined on $T^\ast (\pi(V))$ fixing $z^\ast$ and an open neighbourhood $I\subset \tilde I \subset \R^p$ of $\mu^\ast$ such that for all $\mu \in I$ the maps $\pi|_{\Psi(\Gamma_\mu \cap V)}$ and $\pi|_{\Psi(\Lambda_\mu \cap V)}$ are injective. Let $U=\pi(V)$ and denote $\Psi(\Gamma_\mu \cap V)$ and $\Psi(\Lambda_\mu\cap V)$ again by $\Gamma_\mu$ and $\Lambda_\mu$.\\

There exist 1-forms $\alpha_\mu, \beta_\mu\colon U \to T^\ast U$ such that $\alpha_\mu(U)=\Gamma_\mu$ and $\beta_\mu(U)=\Lambda_\mu$.
The 1-forms $\alpha_\mu$ and $\beta_\mu$ are closed since $\Gamma_\mu$ and $\Lambda_\mu$ are Lagrangian submanifolds. Shrinking $U$ to a simply connected domain around $\pi(z^\ast)$ denoted again by $U$ and the manifolds $\Gamma_\mu$ and $\Lambda_\mu$ to their intersections with $\pi^{-1}(U)$, the 1-forms are exact. We find $f_\mu,g_\mu \colon U \to \R$ such that $\d f_\mu = \alpha_\mu$ and $\d g_\mu=\beta_\mu$. Thus, $\alpha_\mu(x) = \beta_\mu(x) \in \Gamma_\mu \cap \Lambda_\mu$ if and only if $\d (f_\mu-g_\mu) |_x =0$.
\qed

\begin{remark}
If we assume that the families of Lagrangian submanifolds are graphical, then we do not need lemma \ref{lem:graphoverzero} to translate graphical Lagrangian intersection problems to gradient zero problems.
\end{remark}

\begin{prop}\label{prop:rightequivalence}
The presented translation procedure in lemma \ref{lem:graphoverzero} and in the linkage theorem \ref{thm:linkageThm} of local Lagrangian intersection problems $\mathfrak L = \big( (\Lambda_\mu)_{\mu \in  I}, (\Gamma_\mu)_{\mu \in  I},(\mu^\ast,z^\ast) \big) $ to families $(h_\mu)_{\mu \in I}$ of local functions determines $(h_\mu)_{\mu \in I}$ up to a $\mu$-dependent family of diffeomorphisms acting from the right and a $\mu$-dependent family of affine translations acting from the left.
\end{prop}

\textit{Proof.}
As we have already remarked (remark \ref{rem:zerosecinv}), each $\Lambda_\mu$ can be mapped to the zero section by a symplectomorphism which is uniquely determined up to a symplectomorphism leaving the zero-section invariant.
Moreover, notice that the Lagrangian subspaces defined by the normal forms $A_0$ and $B_0$ in lemma \ref{lem:graphoverzero} only depend on the choice of coordinates on the base space $N$.
We achieve $g_\mu=0$ for all $\mu \in I$ in the proof of the Singularity-Lagrangian linkage theorem \ref{thm:linkageThm}. 
A symplectomorphism $\Psi_\mu$ on $T^\ast N$ leaving the zero-section invariant is the cotangent lifted action of a diffeomorphism $\chi_\mu$ on $N$.
It maps the cotangent vectors $\d f_\mu|_{x}$ to $\d (f_\mu \circ \chi_\mu^{-1})|_{\chi_\mu(x)}$. This proves the claim.
\qed\\



An equivalence relation that is coarser than our \textit{fibre-wise right diffeomorphic relation} is \textit{right-left equivalence}. This relation is typically used in singularity theory.
In addition to fibre-wise diffeomorphic changes of variables on the domain, the target space $\R$ of the family of scalar valued maps and the parameters $\mu$ are allowed to transform as well.
See e.g.\ \cite[p.68ff]{lu1976singularity} for definitions and discussions. 
In particular, right-left equivalence allows one to swap $(\Lambda_\mu)_{\mu \in I}$ and $(\Gamma_\mu)_{\mu \in I}$. The rest of the paper refers to the right-left category unless otherwise stated. However, all arguments also work for coarser equivalence relations like contact-equivalence, for instance.\\

Catastrophe theory considers bifurcation behaviour which persists under small perturbations (up to right-left equivalence).
In order to apply classification results from catastrophe theory we
introduce the concept of \textit{genericity}. Roughly speaking, a parametrised family of maps is generic if the family covers all small perturbations.  Introducing a new unfolding parameter does not destroy or lead to new effects. For example, a generic family does not possess any symmetries affecting the bifurcation behaviour which can be destroyed by a small symmetry-breaking perturbation.
This can be made precise as follows.

\begin{definition}[generic family of smooth maps]\label{def:genericfamily}
Let $I \subset \R^p$ and $U \subset \R^n$ be open sets.
A smooth family $(h_\mu)_{\mu \in I}$ of smooth maps $h_\mu \colon U \to \R$ is called a generic family if and only if for each $(\mu^\ast, p^\ast) \in I\times U $ the family $(h_\mu)_{\mu \in I}$ is a versal unfolding of the map germ $h_{\mu^\ast}\colon (U,p^\ast) \to (\R,h_{\mu^\ast}(p^\ast))$.
\end{definition}

Versality of a family in definition \ref{def:genericfamily} is to be understood in the sense of singularity theory \cite[p.72]{lu1976singularity} and means that the family covers all small perturbations.
This means the bifurcation diagrams of the perturbed gradient-zero problems qualitatively look the same as the bifurcation diagram of the unperturbed problem. Versality can be translated into algebraic conditions.\footnote{
Consider the sheaf $\mathcal C^\infty(U)$ of smooth, real-valued maps on a manifold $U$. For $p^\ast \in U$ the elements of the stalk $\mathcal R$ at $p^\ast$ are referred to as map germs and denoted by $g \colon (U,p^\ast) \to (\R,g(p^\ast))$. A smooth family of map germs $(g_\mu)_{\mu \in I} \subset \mathcal R$ with $g_{\mu^\ast}=g$ is called an \textit{unfolding of $g$}, where $\mu^\ast \in I$ is a parameter value and $I$ the parameter space.
The $\R$-algebra structure of $\R$ induces a ring structure on $\mathcal R$.
Let $\mathfrak m$ be the (unique maximal) ideal in $\mathcal R$ of map germs vanishing at $p^\ast$.
For $U = \R$ we denote the corresponding ideal by $\mathfrak m_\R$.
If for a representative $\tilde g$ of the map germ $g$ the differential $\d \tilde g$ is not trivial at $p^\ast$ then any unfolding of $g$ is versal.
In the right-left category it is sufficient to restrict to $g(p^\ast) = 0$.
If $\d \tilde g$ vanishes then $g \in \mathfrak m^2 = \mathfrak m \cdot \mathfrak m$.
The ideal $\mathfrak m_\R$ can be pulled back to $g^\ast \mathfrak m \subset \mathfrak m$ using $g$: if $f \in \mathfrak m_\R$ then $g^\ast(f)=f \circ g$.
For local coordinates $x^i$ on $U$  around $p^\ast$ consider the ideal $\mathcal I_g$ spanned by the map germs of $\frac{\p \tilde g}{\p x^i}$ and by $g^\ast \mathfrak m $.
Define the quotient $Q_g= \mathfrak m/ \mathcal I_g$.
Consider the unfolding 
$g_\mu = g+\sum_{j=1}^k \mu_j g_j$
of $g$ with $g_j \in \mathfrak m$.
Denote the projection of $g_j$ to the quotient $Q_g=\mathfrak m/ \mathcal I_g$ by $[g_j]$. If $[g_1], \ldots, [g_k]$ is a system which $\R$-spans $Q_g$ then the unfolding $(g_\mu)_\mu$ is versal.
}

Thus, the nonvanishing of certain algebraic functions of the derivatives of a function at a point imply that the family is generic in a neighbourhood of that point.
Hence it can, in principle, be checked for any given family.

\begin{definition}[generic Lagrangian intersection problem] Consider a local Lagrangian intersection problem $\big( (\Lambda_\mu)_{\mu \in  I}, (\Gamma_\mu)_{\mu \in  I},(\mu^\ast,z^\ast) \big)$. The problem is referred to as \textit{generic} if an application of the linkage theorem \ref{thm:linkageThm} defines a generic family of maps $(h_\mu)_{\mu}$.
\end{definition}

Now we can formulate some implications of the Singularity-Lagrangian linkage theorem \ref{thm:linkageThm}.

\begin{corollary}\label{cor:genint}
Generic local Lagrangian intersection problems are governed by catastrophe theory. 
\end{corollary}


Moreover, the proof of the Singularity-Lagrangian linkage theorem \ref{thm:linkageThm} tells us under which conditions we leave the setting of catastrophe theory.

\begin{prop}\label{prop:arint}
Consider an intersection problem $\Lambda_\mu \cap \Gamma_\mu$ of two graphical families of submanifolds $(\Lambda_\mu)_\mu\subset T^\ast N$ and $(\Gamma_\mu)_\mu \subset T^\ast N$ where at least one of the families varies through arbitrary, i.e.\ not necessarily Lagrangian, submanifolds of dimension $n$, where $2n$ is the dimension of the ambient space. The intersection problem corresponds to the problem $\alpha_\mu =0$ with a family of 1-forms $\alpha$ on $N$. This is equivalent to the equilibria-of-vector-fields problem/ zeros-of-a-function problem.
\end{prop}

\begin{remark}
References for the zeros-of-a-function problem $F_\mu=0$ for a family $F_\mu \colon \R^m \to \R^k$ are \cite{TopoStability,walllectures}.
\end{remark}

\textit{Proof of proposition \ref{prop:arint}.} The intersection problem corresponds to the problem $\alpha_\mu=0$ for a family of 1-forms $(\alpha_\mu)_\mu$ since the families of submanifolds are graphical. Now we construct an example showing that any family of 1-forms $(\alpha_\mu)_\mu$ can occur even if one of the families, say $(\Gamma_\mu)_\mu$, is special (e.g.\ Lagrangian).
Since the manifolds $\Gamma_\mu$ are graphical, there exist 1-forms $\beta_\mu$ such that $\beta_\mu(N)=\Gamma_\mu$. Set $\Lambda_\mu=(\alpha_\mu+\beta_\mu)(N)$. Now $z \in \Gamma_\mu \cap \Lambda_\mu$ if and only if $\alpha_\mu(\pi(z))=0$.
\qed

\begin{prop}\label{prop:Laggenfix}
A degenerate subcase of proposition \ref{prop:arint} is that the family $\Gamma_\mu$ is Lagrangian and generic while $\Lambda_\mu$ is constant and not Lagrangian. This corresponds to the problem $\d f_\mu = \beta$, where $\beta$ is not closed. 
Although the parameter $\mu$ enters in the same way as in the \textit{gradient-zero-problem}, the occurring bifurcations are \textit{not} generally governed by catastrophe theory.
\end{prop}

\textit{Proof.} If $f_\mu$ is the truncated miniversal unfolding of the hyperbolic umbilic singularity $D_4^+$, i.e.\
\[
f_\mu(x,y)
= x^3+xy^2+\mu_3 (x^2+y^2) + \mu_2 y + \mu_1 x
\]
then the bifurcation behaviour of $\d f_\mu (x,y) = \e \cdot x \d y$ changes qualitatively if $\e >0$, i.e.\ the bifurcation diagram of the perturbed system is not fibrewise diffeomorphic to the diagram of the miniversal unfolding of $D_4^+$. This can be deduced from the fact that only for $\e =0$ there exists a path in the $\mu$-space such that four solutions merge to one solution while in the perturbed system this bifurcation breaks up into two fold bifurcations. This means that $D_4$ together with its unfolding is not a normal form in this problem class, so catastrophe theory does not cover this case.
\qed\\


We can formulate a reverse direction of the Singularity-Lagrangian linkage theorem \ref{thm:linkageThm}:

\begin{prop}\label{cor:linkagereverse}
For each gradient-zero-problem there exists a local Lagrangian intersection problem with the same catastrophe set.
\end{prop}

\textit{Proof.} Let $I = \R^p$ and $f_\mu\colon \R^n \to \R$ be a smooth family of smooth maps.
Let $\Gamma_\mu$ be the image of the 1-form $\d f_\mu$ and $\Lambda_\mu$ be the zero section in $T^\ast \R^n$. The intersection $\{(\mu,z) \in \Gamma_\mu \cap \Lambda_\mu\}$ is a catastrophe set of the family $f_\mu$.\qed

\section{Boundary value problems for symplectic maps}\label{sec:bd-problemssymplmaps}

The following analysis shows how the Singularity-Lagrangian linkage theorem \ref{thm:linkageThm} helps to understand bifurcation phenomena in boundary value problems for symplectic maps.

\subsection{Application of the Singularity-Lagrangian linkage theorem}\label{subsec:appbdsympl}

We will show how to rephrase boundary value problems for symplectic maps as intersection problems of Lagrangian manifolds.

\subsubsection{Lagrangian embedding of symplectic maps}

Let $\phi_\mu \colon (\R^{2n},\omega) \to (\R^{2n},\omega)$ be a family of symplectic maps. We equip $(\R^{2n},\omega)$ with symplectic coordinates $x^1,\ldots,x^n,y_1,\ldots,y_n$ such that $\omega= \sum_{j=1}^n \d x^j \wedge \d y_j$. Then we equip the cotangent bundle $T^\ast \R^{2n}$ with the standard symplectic structure $\d \lambda$ induced by its canonical Liouvillian $\lambda$ defined by \eqref{eq:Liouvillian} with $N=\R^{2n}$. Darboux coordinates are given by $q^1,\ldots,q^{2n},p_1,\ldots,p_{2n}$ with
\begin{align*}
q^i \colon T^\ast N \to \R \quad &q^i = x^i \circ \pi\\
q^{i+n} \colon T^\ast N \to \R \quad &q^{i+n}=y_i \circ \pi\\
p_i \colon T^\ast N \to \R \quad & \beta \mapsto \beta\left(\left.\frac{\p}{\p x^i}\right|_{\pi(\beta)}\right)\\
p_{i+n} \colon T^\ast N \to \R \quad & \beta \mapsto \beta\left(\left.\frac{\p}{\p y_i}\right|_{\pi(\beta)}\right)
\end{align*}
for $i \in \{1,\ldots,n\}$.
The canonical Liouvillian $\lambda$ has the coordinate expression
\[ 
\lambda = \sum_{i=1}^{2n} p_i \d q^i = \sum_{i=1}^{n} \left( p_i \d q^i + p_{n+i} \d q^{n+i}\right).
\]
Denote the composition of $\phi_\mu$ with a projection to the $x$-coordinates by $\phi_\mu^X$  and to the $y$-coordinates to $\phi_\mu^Y$.
The immersion $\iota \colon \R^{2n} \to T^\ast \R^{2n}$ with 
\begin{equation}\label{def:iota}
\iota_\mu(x,y) 
= \left(x, \phi_\mu^Y (x,y),y, \phi_\mu^X(x,y) \right)
\end{equation}
is Lagrangian, i.e.\ $\d \iota^\ast\lambda =0$, because $\phi_\mu$ is symplectic. Define $\Gamma_\mu = \iota_\mu(\R^{2n})$.

\subsubsection{Embedding of boundary conditions}
Families of boundary conditions for $\phi_\mu$ can be represented as families of submanifolds $\Lambda_\mu \subset T^\ast\R^{2n}$ such that solutions correspond to the elements in the intersections $\Gamma_\mu \cap\Lambda_\mu $. For well-posed\footnote{The degrees of freedom of the problem equals the number of independent boundary conditions.} boundary value problems, the manifolds $\Lambda_\mu$ are $2n$-dimensional.

\begin{definition}[Lagrangian boundary conditions]
Boundary conditions for which the corresponding submanifolds $\Lambda_\mu \subset T^\ast\R^{2n}$ are Lagrangian are called \textit{Lagrangian boundary conditions}.
\end{definition}

Examples are

\begin{itemize}
\item
the boundary condition
\[\Lambda_\mu = \{ (x^\ast,Y,y,X^\ast)|y,Y \in \R^n\}\]
for fixed $x^\ast,X^\ast \in \R^n$ representing the two point-boundary value problem for $\phi_\mu$ in the $x$-coordinates, i.e.\
\[
x = x^\ast, \quad \phi^X_\mu(x,y) = X^\ast.
\]

\item
Choosing $\Lambda_\mu$ to be the zero section corresponds to the problem
\[ y=0, \quad \phi^X_\mu(x,y) = 0.\]

\item
Periodic boundary conditions $\phi_\mu(x,y)=(x,y)$ are represented by
\[\Lambda_\mu = \{ (x,y,y,x)|x,y \in \R^n\}.\]
\end{itemize}

An example for a \textit{non}-Lagrangian boundary condition is $\phi_\mu(x,y)=(y,x)$ with
\[\Lambda_\mu = \{ (x,x,y,y)|x,y \in \R^n\}.\]


The representation of boundary conditions as submanifolds in the cotangent bundle depends, of cause, on the specific choice of the embedding of the symplectic maps. The above examples are given with respect to the choice $\iota$ as defined in \eqref{def:iota}. However, whether a boundary condition is Lagrangian or not does not depend on the specific map $\iota$ as long as we require $\iota$ to be a Lagrangian embedding. \\

The requirement that the family of symplectic maps is globally defined on $N=\R^{2n}$ was imposed for convenience.
Indeed, the discussion applies to arbitrary families of symplectic maps $(\phi_\mu \colon X \to X')_{\mu \in I}$: since we describe \textit{local} behaviour around some $z^\ast \in X$ and $\mu^\ast \in I$ we can assume that $X$ and $X'$ admit Darboux coordinates centred at $z^\ast$ and $\phi_{\mu^\ast}(z^\ast)$ respectively. The discussion then applies to the family $(\phi_\mu)_{ \mu \in I}$ expressed in centred Darboux coordinates. If we want to consider the local behaviour of a smooth family $(\phi_\mu)_{\mu \in \tilde I}$ of symplectic maps on a symplectic manifold $N$ and a parameter space $\tilde I$ locally around $\mu^\ast \in \tilde I$ and $z^\ast \in N$, for instance, then we take $I$ as such a small neighbourhood of $\tilde I$ around $\mu^\ast \in \tilde I$ and $X$ as such a small neighbourhood of $z^\ast \in N$ that $X$ and $X' =\bigcup_{\mu \in I}\phi_\mu (X)\subset N$ both admit Darboux coordinates. 

\subsubsection{Towards a coordinate-free description of embedding of symplectic maps}\label{subsub:coordfree}
In \eqref{def:iota} we have defined an embedding $\iota$ of a family of symplectic maps in Darboux coordinates because boundary conditions are typically given in coordinates. However, the following viewpoint is natural in mechanical Hamiltonian systems and it is convenient when describing the interaction of generating functions and symmetries.\\

Let $U$ be a finite-dimensional vector space and $U^\ast$ its dual vector space. On the direct sum $V = U \oplus U^\ast$ we consider the symplectic form
\[
\omega\big( (u,u^\ast),(v,v^\ast) \big) = u^\ast(v)-v^\ast(u).
\]
Denote a copy of the symplectic vector space $(V,\omega)$ by $(\underline V = \underline U \oplus \underline U, \underline \omega)$. Consider the symplectic space $\mathcal U = \big(V \oplus \underline V, \omega \oplus (- \underline \omega)\big)$ with symplectic form
\[
\omega \oplus (- \underline \omega) = \mathcal P^\ast \omega - \underline {\mathcal P}^\ast \underline \omega,
\]
where $\mathcal P \colon \mathcal U\to V$ and $\underline{\mathcal P} \colon \mathcal U\to \underline V$ denote the projections.

\begin{lemma}
The symplectic space $\mathcal U$ can be identified canonically\footnote{in a coordinate invariant way} and symplectomorphically with the cotangent bundle $T^\ast(U \oplus \underline U^\ast)$ equipped with the canonical symplectic form $\d \lambda$, where $\lambda$ is the Liouvillian-1-form (defined as in \eqref{eq:Liouvillian}).
\end{lemma}

\textit{Proof.} The identification is obtained in two steps:

\begin{itemize}
\item We can identify $T^\ast(U \oplus \underline U^\ast) \cong T^\ast U \oplus T^\ast \underline U^\ast$. Consider the natural projections
\[
\pi \colon U\oplus \underline U^\ast \to U, \quad 
\underline \pi \colon U\oplus \underline U^\ast \to \underline U^\ast.
\]
The identification $T^\ast U \oplus T^\ast \underline U^\ast \cong T^\ast(U \oplus \underline U^\ast)$ is given by the isomorphism
\[
(\alpha,\beta) \mapsto \pi^\ast \alpha 
+\underline \pi^\ast  \beta 
\]

\item
We can identify $T^\ast U \oplus T^\ast \underline U^\ast \cong V \oplus \underline V$.
Consider the bundle projection $\kappa \colon T^\ast U \to U$. For each $u \in U$ we denote the canonical, linear identification  of the vector space $U$ and its tangent space $T_uU$ by $L_u\colon U \to T_uU$. The identification $T^\ast U \cong U \oplus U^\ast$ is obtained as
\[
\alpha \mapsto (\kappa(\alpha), \alpha \circ L_{\kappa(\alpha)})
\]

The identification $T^\ast \underline U^\ast \cong \underline U^\ast \oplus \underline {(U^\ast)^\ast}$ is obtained in an analogous way. Using the linear, canonical identification of $U$ and its bi-dual space $\mathrm{eval}\colon U \to (U^\ast)^\ast$ we obtain $T^\ast \underline U^\ast \cong \underline U^\ast \oplus \underline U$.
\end{itemize}

The fact that the described canonical identification $T^\ast(U \oplus \underline U^\ast)$ is symplectic is most easily seen by expressing the maps in symplectic coordinates in which the identifications appear trivial.
\qed\\

Let $B$ and $\underline B$ denote open neighbourhoods of $0 \in V$ and $0 \in \underline V$. The graph 
\[
\Gamma = \{(v,\phi(v)) | v \in V \cap B\} \subset \mathcal U.
\]
of a symplectic map $\phi \colon (V \cap B,\omega) \to (\underline V \cap \underline B,\underline \omega)$ is a Lagrangian submanifold. It gives rise, therefore, to a Lagrangian submanifold in $T^\ast(U \oplus \underline U^\ast)$. Thus, we obtain a coordinate free description to embed the graph of a local symplectic map on $V=U \oplus \underline U^\ast$ into $T^\ast(U \oplus \underline U^\ast)$. In case $\Gamma$ is graphical, there exists $S \colon U \oplus \underline U^\ast \to \R$ such that $\Gamma$ is the image of $\d S$. We may use $S$ as a generating function.\\

In applications, the splitting $U \oplus \underline U^\ast$ occurs naturally: in Hamiltonian mechanics, $U$ models a configuration space and $\underline U^\ast$ the space of momenta. However, in the more general setting of a symplectic map $\tilde \phi$ defined on a symplectic manifold $M$ locally around $z \in M$, we would use Darboux patches around $z$ and $\tilde \phi(z)$ to obtain a symplectic map $\phi \colon (V \cap B,\omega) \to (\underline V \cap \underline B,\underline \omega)$ that fits into the analysed setting. This, of cause, involves choices.

\subsubsection{Summary of the translation procedure.}
We summarise the key ideas of the translation of parameter dependent Lagrangian boundary value problems into gradient-zero-problems:

\begin{itemize}
\item
Symplectic maps and Lagrangian boundary conditions constitute Lagrangian submanifolds in a cotangent bundle.
\item
After shrinking the involved manifolds and the parameter space and after applying a symplectic transformation all involved Lagrangian submanifolds can assumed to be graphical.
\item
The Lagrangian submanifolds can be written as images of exact 1-forms and the primitives of the forms coming from the boundary conditions can be substracted from the primitives related to the symplectic maps.
\item
We obtain a family of smooth functions whose critical points correspond to solutions of the boundary value problem.
\end{itemize}

(Notice that the motivational example from section \ref{subsec:motivation} shows a slightly different viewpoint where the 1-forms are defined on the Lagrangian graph $\Gamma_\mu \subset \R^2 \times \R^2$ and the primitive $\alpha$ of the symplectic form $\omega_\times$ is chosen in view of the Dirichlet boundary value problem. )

\subsubsection{Translation of results from catastrophe theory}

\begin{definition}[generic family of symplectic maps]\label{def:gensympl}
A family $(\phi_\mu)_{\mu \in I}$ of symplectic maps $\phi_\mu\colon (X,\omega) \to (X',\omega')$ is referred to as generic if for each $(\mu^\ast,x^\ast)\in I \times X$ there exist an open neighbourhood $\tilde I \times \tilde X \subset I \times X$ such that there exists a generic family of maps $(h_\mu)_{\mu \in \tilde I}$ (see definition \ref{def:genericfamily}) generating\footnote{In the sense of generating functions of symplectic maps, see e.g.\ \cite{daSilva2008}} the symplectic maps $(\phi_\mu|_{\tilde X})_{\mu \in \tilde I}$.
\end{definition}

Now we can use corollary \ref{cor:genint}, proposition \ref{prop:arint} and proposition \ref{prop:Laggenfix} to obtain the following statements.

\begin{prop}\label{prop:linkCata} Let $\phi_\mu\colon (X,\omega) \to (X',\omega')$ be a generic family of symplectic maps with Lagrangian embeddings $\Gamma_\mu=\iota_\mu(X)$ defined by \eqref{def:iota} using locally centred Darboux coordinates. Consider a boundary value problem for $\phi_\mu$ where the boundary conditions are represented by a family of submanifolds $\Lambda_\mu$ in a neighbourhood of $0 \in T^\ast\R^{2n}$.
\begin{itemize}
\item
If the family of submanifolds $(\Lambda_\mu)_\mu$ is constant and Lagrangian, then the intersection $\{(\mu,z) \in \Gamma_\mu \cap \Lambda_\mu\}$ corresponds to a catastrophe set in catastrophe theory.

\item
If the family $(\Lambda_\mu)_\mu$ varies through arbitrary $2n$-dimensional submanifolds of a neighbourhood of $0 \in$ $T^\ast\R^{2n}$ then the intersection $\{(\mu,z) \in \Gamma_\mu \cap \Lambda_\mu\}$ corresponds to the zeros-of-a-function problem / equilibria-of-vector-fields problem.

\item
If the family$(\Lambda_\mu)_\mu$ is a constant non-Lagrangian submanifold of dimension $2n$ then the intersection $\{(\mu,z) \in \Gamma_\mu \cap \Lambda_\mu\}$ is, generally speaking, not governed by catastrophe theory.
\end{itemize}

\end{prop}

The motivational example (figure \ref{fig:cusp-example-problem}) is an instance of a constant Lagrangian boundary condition (Dirichlet conditions) and a generic family of symplectic maps which is given as a family of Hamiltonian diffeomorphisms.

\subsection{Obstructions for bifurcations in low dimensions}\label{subsec:nogobifur}

%

In the motivational example we considered a Dirichlet problem for the flow map of a planar Hamiltonian system (section \ref{subsec:motivation}) and found a cusp bifurcation (figure \ref{fig:cusp-example-problem}). By the (real) ADE-series classification in the right-left category obtained by Arnold \cite[p.33]{Arnold1}, the cusp bifurcation belongs to the A-series. 
The general classification provides some motivation to introduce more parameters in the planar Hamiltonian and chase up more complicated bifurcations. Clearly, we cannot find bifurcations whose germs have normal forms with more than two variables because the number of variables used for the normal forms in the ADE classification tables is minimal. However, the normal forms of $D$ or $E$ series bifurcations only use two variables. Do they occur in planar Dirichlet problems?\\

In the following we will prove that, no matter how many parameters we introduce, only $A$ series bifurcations are possible in the planar Dirichlet problem. Indeed, we would need an at least four-dimensional system to find $D$ or $E$ series bifurcations. This is also true for Dirichlet problems for arbitrary symplectic maps which do not necessarily arise as Hamiltonian flows. The reason why we have to double the dimension can be seen from the geometric picture that we develop in this section.\\

The geometric idea is the following: by proposition \ref{cor:linkagereverse} each type of local Lagrangian intersection can be achieved with two families of Lagrangian submanifolds $(\Gamma_\mu)_{\mu \in I}, (\Lambda_\mu)_{\mu \in I} \subset T^\ast \R^{2n}$. However, if the boundary condition $\Lambda_\mu$ is constant in $\mu$ and if  $(\Gamma_\mu)_{\mu \in I}$ is a Lagrangian embedding of the graphs of symplectic maps (we may take \eqref{def:iota}) then this restricts the way $\Gamma_\mu$ and $\Lambda_\mu$ can intersect. Depending on the position of $\Lambda_\mu$ in $T^\ast \R^{2n}$, the manifolds are not allowed to touch in a way that the intersection of their tangent spaces at the touch point is of maximal dimension. This prohibits certain bifurcations.\\

The following proposition follows from the proof of the Singularity-Lagrangian linkage theorem \ref{thm:linkageThm}.

\begin{prop}\label{prop:obstr}
Let $\big( (\Lambda_\mu)_{\mu \in  I}, (\Gamma_\mu)_{\mu \in  I},(\mu^\ast,z^\ast) \big)$ be a local Lagrangian intersection problem.
If $m$ is the dimension of the intersection of the tangent spaces $T_{z^\ast}\Lambda_{\mu^\ast}$ and $T_{z^\ast}\Gamma_{\mu^\ast}$
then only those singularities can occur
at $(\mu^\ast,z^\ast)$
which can be obtained in the gradient-zero-problem $\nabla g_\mu=0$ for a family of smooth maps $g_\mu \colon \R^k \to \R$ with $k\le m$.
\end{prop}

\begin{remark}
If $m=0$ then the above proposition states that there is no singularity at $(\mu^\ast,z^\ast)$.
\end{remark}

Let us prove statements for Dirichlet problems and for periodic boundary conditions. 

\begin{prop}\label{prop:dimrestr}
Let $x^\ast, X^\ast \in \R^{2n}$. In the Dirichlet problem
\begin{equation*}
x=x^\ast, \quad \phi^X_\mu(x,y)=X^\ast
\end{equation*}
for a smooth family of symplectic maps $\phi_\mu \colon \R^{2n} \to \R^{2n}$ only those singularities occur which can be obtained in the gradient-zero-problem $\nabla g_\mu=0$ for a smooth family of smooth maps $g_\mu \colon \R^k \to \R$ with $k \le n$.
\end{prop}

In a generic setting, proposition \ref{prop:dimrestr} allows us to use the ADE-classification \cite[p.33]{Arnold1} to obtain statements about which singularities can occur.
We can conclude

\begin{corollary}\label{cor:OnlyAseries}
In Dirichlet problems for generic families of symplectic maps on $T^\ast \R$ only A-series bifurcations occur.
\end{corollary}

\begin{corollary}\label{cor:OnlyModality0}
In Dirichlet problems for generic families of symplectic maps on $T^\ast \R^2$ with at most 7 parameters only singularities of modality 0 (simple singularities) occur.
\end{corollary}

\textit{Proof of proposition \ref{prop:dimrestr}.}
Since symplectic maps are local diffeomorphisms, $\d \iota_\mu$ is a bundle isomorphism between $T\R^{2n}$ and $T \Gamma_\mu$ which in the frame $\frac{\p}{\p q^1},\ldots,\frac{\p}{\p q^{2n}},\frac{\p}{\p p_1},\ldots,\frac{\p}{\p p_{2n}}$ is given by
\begin{equation}\label{eq:Jacobiiota}
D \iota_\mu = \begin{pmatrix}
\Id_{ n}&0\\ 
D_x \phi^Y_\mu & \D_y \phi^Y_\mu\\
0 & \Id_{n} \\ 
D_x \phi^X_\mu & \D_y \phi^X_\mu\\
\end{pmatrix}.
\end{equation}
Consider the Dirichlet problem $\Lambda_\mu = \{ (x^\ast,Y,y,Y^\ast)|y,Y \in \R^n\}$ whose tangent spaces are spanned by the columns of
\[
\begin{pmatrix}
0&0\\
\Id_{n}&0\\
0&\Id_{n }\\
0&0
\end{pmatrix}.
\]
Therefore, at an intersection of $\Gamma_\mu$ and $\Lambda_\mu$ the tangent spaces intersect in an isotropic linear space with dimension at most $n$. 
We conclude that only those singularities can occur which in catastrophe theory admit a description with at most $n$ variables.
\qed\\



Periodic boundary conditions, on the other hand, do not impose restrictions, as the following proposition shows.

\begin{prop}\label{prop:allbifurper} In periodic boundary value problems $\phi_\mu(x,y)=(x,y)$ for smooth families of symplectic maps in $2n$ variables all singularity bifurcations which admit descriptions in at most $2n$ variables occur.
\end{prop}

\textit{Proof.} Consider the vector space $\R^{2n}$ with coordinates $x^1,\ldots,x^n,Y_1,\ldots,Y_n$.
Consider a family of scalar valued maps $\tilde h_\mu(t_1,\ldots,t_r)$ in $r \le 2n$ variables which is fully reduced at $(\mu,t)=(0,0)$, i.e.\ $\tilde h_0(0)=0$, $\nabla \tilde h_0 (0)=0$ and vanishing Hessian matrix $\mathrm{Hess}\, \tilde h_0 (0)=0$. 
The map
\[
h_\mu(t_1,\ldots,t_{2n}) = \tilde h_\mu(t_1,\ldots,t_r) + t^2_{r+1}+\ldots+t^2_{2n}
\]
is stably equivalent / reducible to $\tilde h_\mu$ (see \cite[p.12]{Arnold1},\cite[p.86]{lu1976singularity} for definitions).
The matrix $\frac{\p^2 g_\mu}{\p x \p Y}$ of mixed derivatives of the map
\[
g_\mu(x,Y)=h_\mu(x,Y)+\sum_{j=1}^n x^j Y_j
\]
is the identity matrix at $(\mu,x,Y)=(0,0,0)$. Now we use $g_\mu$ as a generating functions to obtain symplectic maps $\phi_\mu$: by the implicit function theorem, the system of equations $y = \nabla_x g_\mu(x,Y)$ is solvable for $Y$ defining a smooth function $Y_\mu(x,y)$ with $Y_\mu(0,0)=0$ near $(\mu,x,Y)=(0,0,0)$. Defining $X_\mu(x,y)=\nabla_Y g (x,Y_\mu(x,y))$ we obtain a map $\phi_\mu (x,y)= (X_\mu(x,y),Y_\mu(x,y))$ locally around $0 \in \R^{2n}$. The map $\phi_\mu$ is symplectic with respect to $\omega = \sum_{j=1}^n \d x^j \wedge \d y_j$.
Locally around $(\mu,x,Y)=(0,0,0)$ the image of $\d g_\mu \colon \R^{2n} \to T^\ast\R^{2n}$ coincides with the image of the Lagrangian immersion $\iota_\mu$ defined in \eqref{def:iota} w.r.t.\ the symplectic map $\phi_\mu$. \\

Defining $f(x,Y)=\sum_{j=1}^n x^j Y_j$ the image $\Lambda$ of $\d f \colon \R^{2n} \to T^\ast\R^{2n}$ corresponds to periodic boundary conditions (see section \ref{subsec:appbdsympl}). Since $\d g_\mu - \d f = \d h_\mu$, the periodic boundary value problem for $\phi_\mu$ shows the same bifurcation behaviour as $h_\mu$.
\qed\\

\textit{Example.} Applying the construction in proposition \ref{prop:allbifurper} to the truncated miniversal deformation of the hyperbolic umbilic singularity $D_4^+$ (see table on page \pageref{tabel:elementaryCats}), i.e.\ to
\[
h_\mu(t_1,t_2)=t_1^3+t_1t_2^2+\mu_3(t_1^2-t_2^2)+\mu_2 t_2+\mu_1 t_1,
\]
we obtain the generating function
\[
g_\mu(x,Y)=x^3+xY^2+\mu_3(x^2-Y^2)+\mu_2 Y+\mu_1 x+xY
\]
and the symplectic map $\phi_\mu(x,y)=\big(X_\mu(x,y),Y_\mu(x,y)\big)$ with
\begin{align*}
X_\mu(x,y)&=(x-\mu_3)\left(-1+2\sqrt{-3x^2-2x\mu_3+y-\mu_1+\frac 14}\right)+\mu_2+x\\
Y_\mu(x,y)&= -\frac 12 + \sqrt{-3x^2-2x\mu_3+y-\mu_1+\frac 14}.
\end{align*}
The map $\phi_\mu$ is well-defined near $(\mu_1,\mu_2,\mu_3,x,y)=(0,0,0,0,0)$ and has a fixed point if and only if $\nabla h_\mu(x,y)=0$.\\

If we do not restrict the dimension of the phase space then we can obtain each singularity occurring in the gradient-zero-problem for a family $h_\mu \colon \R^n \to \R$ in any Dirichlet problem
\begin{equation}\label{eq:DirichletinProp}
x=x^\ast, \quad \phi^X_\mu(x,y)=X^\ast
\end{equation}
with $x^\ast, X^\ast \in \R^n$. Indeed, for any type of a Lagrangian boundary value problem and any gradient-zero singularity we can construct a family of symplectic maps defined on a sufficiently high dimensional space such that the boundary value problem undergoes the same bifurcation as the gradient-zero-problem: each Lagrangian boundary condition locally has a generating function. We formulate the following proposition for the generating function suitable for the Dirichlet problem \eqref{eq:DirichletinProp}. 

\begin{prop}\label{prop:allbifur}
Let $h_\mu \colon \R^n \to \R, \, y \mapsto h_\mu(y)$ be a family of smooth maps locally defined around $(\mu,y)=(0,0)$ and let $B_\mu \colon \R^n \times \R^n \to \R, \, (y,Y)\mapsto B_\mu(y,Y)$ be a smooth family locally defined around $(\mu,y,Y)=(0,0,0)$. There exists a family of symplectic maps $\phi_\mu \colon \R^{2n} \to \R^{2n}, \, (x,y)\mapsto \phi_\mu(x,y)$ locally defined around $(\mu,x,y)=(0,0,0)$ such that the boundary value problem
\[
\phi_\mu^X(\nabla_y B_\mu(y,Y),y) = -\nabla_Y B_\mu(y,Y)
\]
%
locally around the origin shows the same bifurcation behaviour as the gradient-zero problem $\nabla h_\mu(y)=0$.
\end{prop}

\begin{remark}
For the Dirichlet problem \eqref{eq:DirichletinProp} we can use
\[
B_\mu(y,Y) = \langle y,x^\ast\rangle - \langle Y,X^\ast\rangle
\]
in the above statement. Here $\langle . , .\rangle$ denotes the scalar product on $\R^n$. 
\end{remark}

\textit{Proof of proposition \ref{prop:allbifur}.}
Define
\[
g_\mu (y,Y) = h_\mu(y)+c \sum_{j=1}^n (Y_j + y_j)^2 + B_\mu(y,Y),
\]
where $c\in \R$ is a constant such that $\frac{\p^2 B_\mu}{\p y \p Y} + 2 c I_n$ is invertible near $(\mu,y,Y)=(0,0,0)$. Since $\frac{\p^2 g_\mu}{\p y \p Y}$ is invertible, the system of equations $x=\nabla_yg_\mu(y,Y)$ locally defines $Y_\mu(x,y)$ by the implicit function theorem. We obtain the required family of symplectic maps as \[\phi_\mu(x,y)=\big(-\nabla_Yg_\mu(y,Y_\mu(x,y)),Y_\mu(x,y)\big).\]
\qed\\

We may formulate a corollary analogous to proposition \ref{cor:linkagereverse} for boundary value problems:

\begin{corollary}\label{cor:summaryLink}
For each (generic) gradient-zero problem there exists a (generic) family of symplectic maps such that the boundary value problem shows the same bifurcation behaviour as the gradient-zero problem.
\end{corollary}

\begin{remark}
The above corollary \ref{cor:summaryLink} holds with or without the parenthesized word ``generic''. Let us emphasise at this point that the translation procedure of Lagrangian boundary value problems to gradient-zero problems (lemma \ref{lem:graphoverzero}, theorem \ref{thm:linkageThm}, proposition \ref{prop:rightequivalence}, proposition \ref{prop:arint} and backwards - proposition \ref{cor:linkagereverse}) is obtained without a restriction to a generic setting. Genericity is needed to link results to the known ADE classification of singularities in catastrophe theory as in corollary \ref{cor:genint} and proposition \ref{prop:linkCata}. However, proposition \ref{prop:obstr}, proposition \ref{prop:dimrestr} and proposition \ref{prop:allbifurper}, dealing with restrictions on which bifurcations can occur, hold for general systems. The implications for a generic settings are formulated in the corollaries \ref{cor:OnlyAseries} and \ref{cor:OnlyModality0}.
\end{remark}

\subsection{Effects of complete integrability}\label{subsec:completeIntegrable}

In the previous section we have seen that Dirichlet problems for symplectic maps allow fewer types of bifurcation then one might naively expect from the ADE-classification results. In this section they will surprise us with more bifurcations than expected.\\

In applications symplectic maps often arise as time-$\tau$-maps of Hamiltonian flows. Those arising in completely integrable 
Hamiltonian systems form an important subclass. Families of Hamiltonian diffeomorphisms in completely integrable systems show more bifurcations than expected from the ADE-classification in certain very common boundary value problems like homogeneous Dirichlet problems. We will present a numerical example first and then develop a general model to explain the \textit{periodic pitchfork bifurcation}, which is caused by complete integrability. It provides a non-trivial example for the effects of symmetries on bifurcation behaviour, which are analysed in \ref{subsec:symeffects}.

\subsubsection{Numerical example.}\label{subsubsec:numex}
\begin{figure}
\begin{center}
\includegraphics[width=0.7\textwidth]{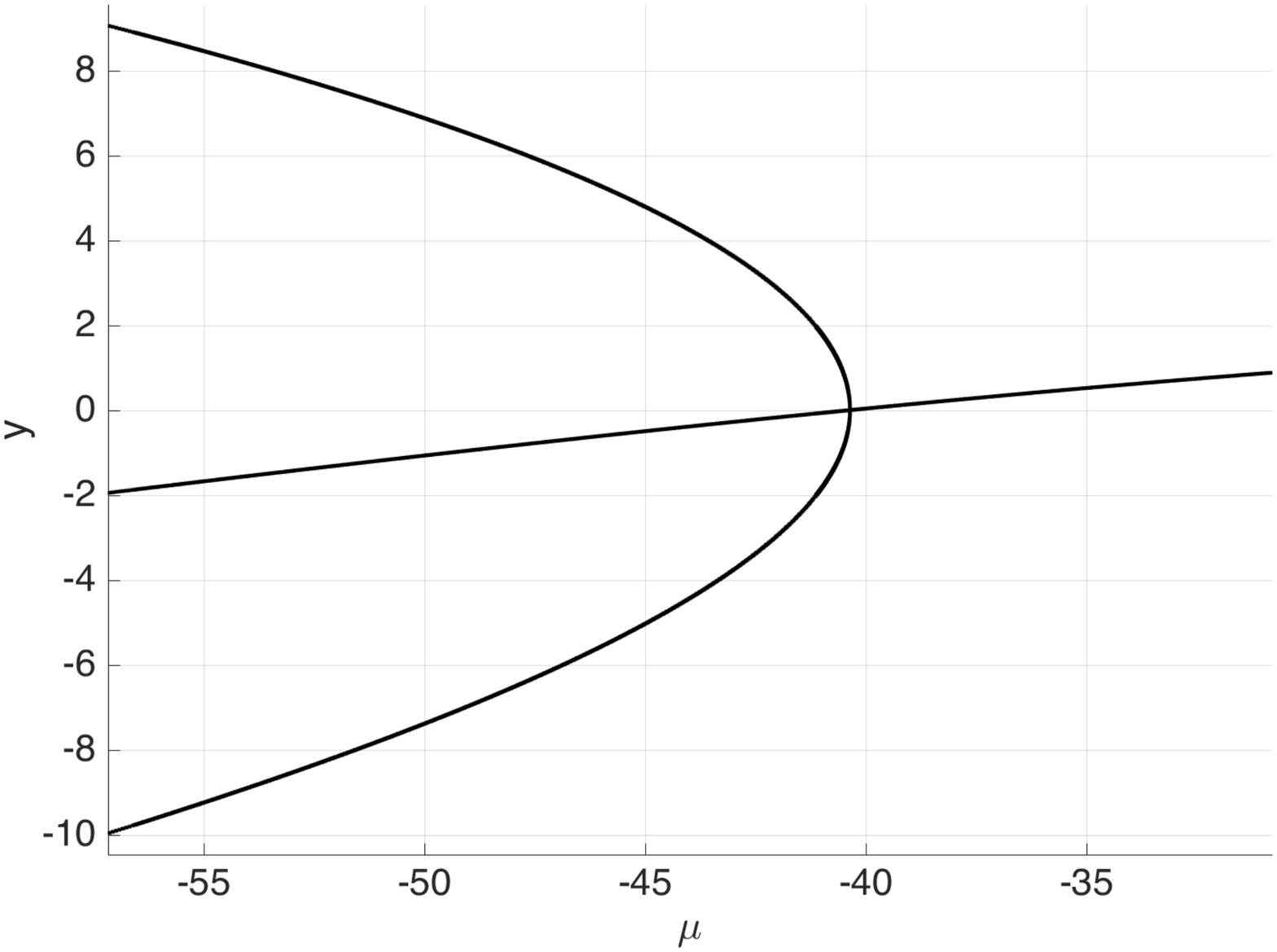}
\caption{The diagram shows a pitchfork bifurcation in the 2-point-boundary value problem $x(0)=1=x(1)$ for the time-1-map of the Hamiltonian $H(x,y) = y^2 +x^3+\mu x + 0.01 y^3$. The flow map is obtained numerically using the (symplectic) leapfrog method. The boundary value problem is solved using a shooting method.}\label{fig:HamPitchInt}
\end{center}
\end{figure}

\begin{figure}
\begin{center}
\includegraphics[width=0.49\textwidth]{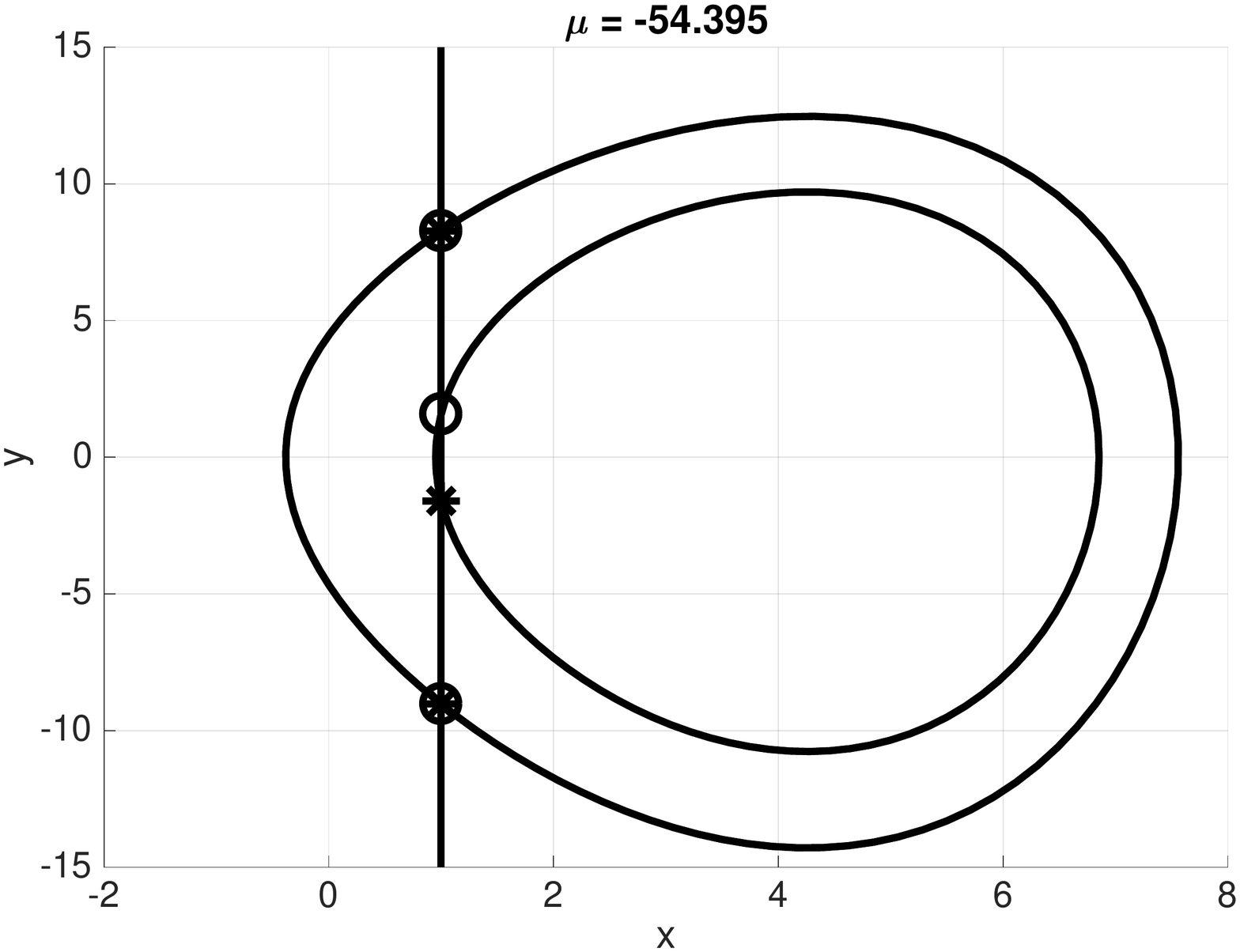}
\includegraphics[width=0.49\textwidth]{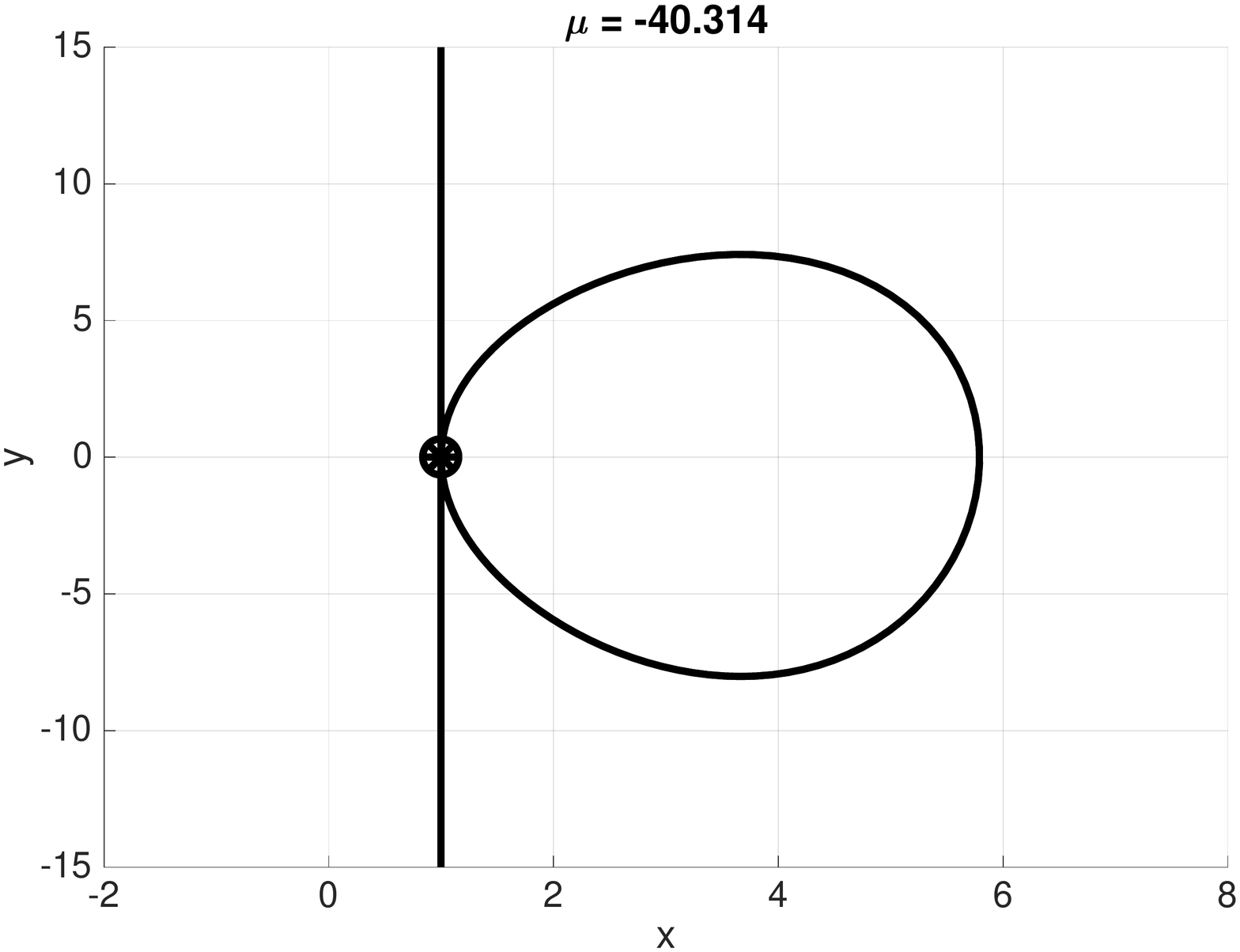}
\includegraphics[width=0.49\textwidth]{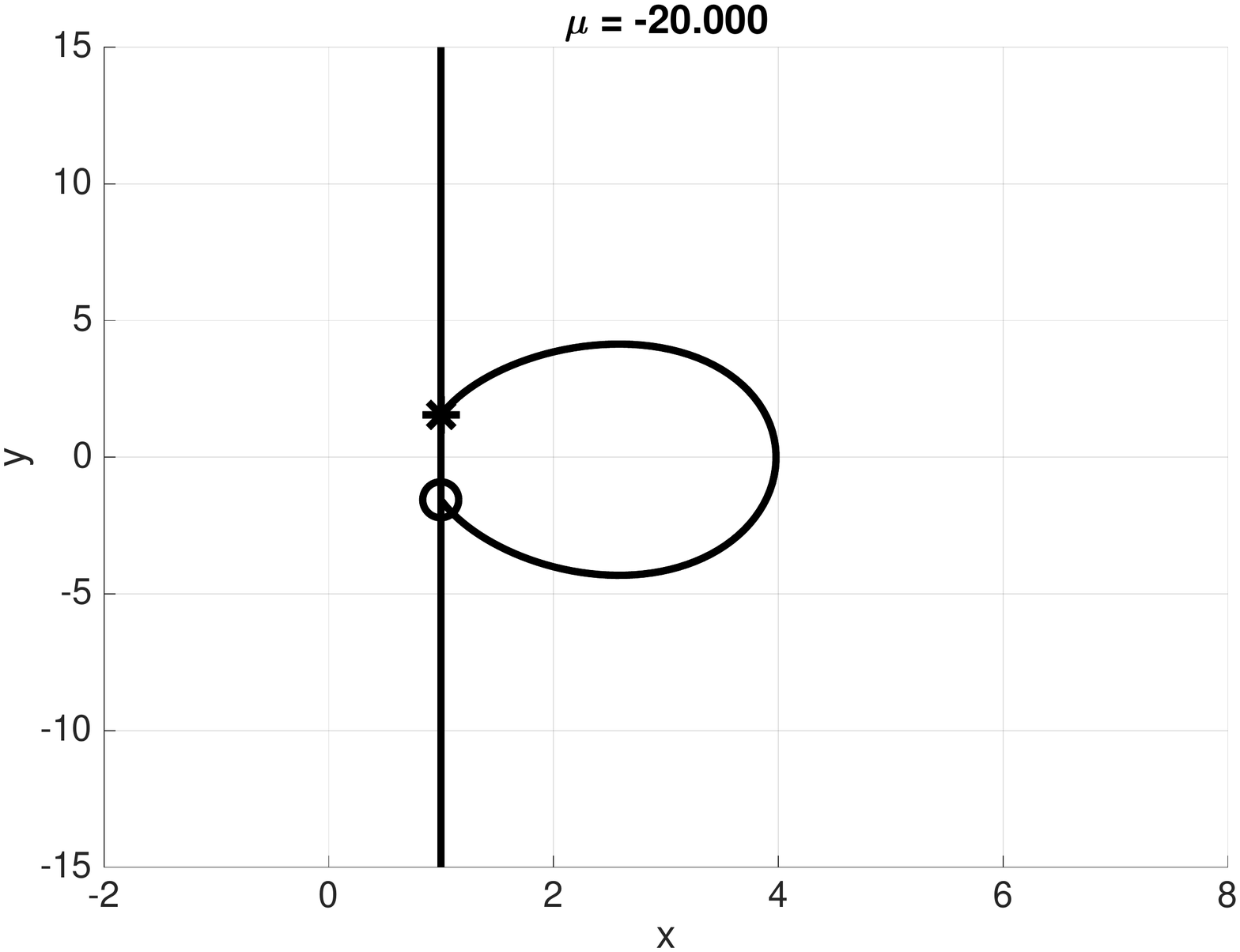}
\caption{The diagram shows the orbits corresponding to the solution branches in figure \ref{fig:HamPitchInt} for different values of $\mu$. It illustrates the mechanism of the periodic pitchfork bifurcation. The orbits are obtained using the (symplectic) leapfrog method. The marker $\ast$ denotes the start point of the orbit and $o$ the end point.}\label{fig:HamPitchIntPhase}
\end{center}
\end{figure}

Planar Hamiltonian systems are completely integrable. 
Let us consider the Hamiltonian\footnote{The term $0.01 y^3$ is included to break global time-reversal symmetry, see also section \ref{subsec:timereversalsym}.}
\[
H(x,y) = y^2 + 0.01 y^3 +x^3+\mu x 
\]
on the symplectic vector space $(\R^2,\d x \wedge \d y)$. Figure \ref{fig:HamPitchInt} shows the bifurcation diagram for the 2-point boundary value problem
\begin{equation}\label{eq:bdProblem}
x(0)=1, \quad x(1)=1
\end{equation}
for the time-1-map of the Hamiltonian flow. We see a pitchfork bifurcation. Figure \ref{fig:HamPitchIntPhase} shows the orbits involved in the bifurcation. The marker $\ast$ is used to indicate the start point of the orbit and $o$ is used for the end point. The black solid line indicates the boundary condition $G =\{1\} \times \R$ in the phase space. At $\mu=-53.306$ there are three solutions to the boundary value problem.
The two solution from the outer branches of the pitchfork in figure \ref{fig:HamPitchInt} come from a periodic orbit with period 1 crossing $G$ twice.
The solution from the inner branch corresponds to another periodic orbit with period slightly smaller than 1.
Increasing $\mu$, all three solutions merge to an orbit of period 1 tangent to $G$.
Increasing $\mu$ further to $-20$, there exists one solution of the boundary value problem which comes from an orbit with period slightly greater than 1.\\

According to the ADE classification in catastrophe theory, the only bifurcation that occurs in generic 1-parameter boundary value problems is the fold bifurcation. This means that the family of Hamiltonian diffeomorphisms is not generic around the pitchfork point $(\mu^\ast,(x^\ast,y^\ast))$ in the sense of definition \ref{def:gensympl}. Indeed, the bifurcation makes use of the fact that periodic orbits occur as 1-parameter families in planar Hamiltonian systems. The trick is, basically, that \textit{all} intersection points of $G$ with a periodic orbit of period $\tau$ give rise to solutions. Let us, therefore, call the phenomenon \textit{periodic pitchfork bifurcation}.
As we will see in the next section, in higher dimensional completely integrable systems this generalises to the fact that all intersections of $G$ with a Liouville torus of period $\tau$ give rise to solutions. 

\subsubsection{Theoretical model for periodic pitchforks in completely integrable systems}\label{subsubsec:periodicpitchforkexpl}

We generalise the observations from the numerical example to higher dimensions and to a more general class of boundary conditions.

\begin{definition}[symmetrically separated Lagrangian boundary conditions]\label{def:symbd}
Consider a family of symplectic maps $\phi_\mu\colon X \to X$ and let $G_\mu$ be a family of Lagrangian submanifolds in $X$. We consider the boundary condition
\[z \in G_\mu \quad \text{and}\quad \phi_\mu(z)\in G_\mu.\]
Let us refer to this type of boundary condition as \textit{symmetrically separated Lagrangian boundary condition}.
\end{definition}

\textbf{Example.} The Dirichlet boundary condition \eqref{eq:bdProblem} used in the numerical example (section \ref{subsubsec:numex}) is symmetrically separated in the sense of definition \ref{def:symbd}. The corresponding manifold $G_\mu$ is given as the line $G=\{1\}\times \R$. More generally, Dirichlet conditions are symmetrically separated boundary conditions if and only if the start- and end-point coincide.

\begin{remark} In contrast to the definition of Lagrangian boundary conditions, the manifold $G_\mu$ is required to be a Lagrangian submanifold of the phase space $X$ (and not of the cotangent bundle $T^\ast X$). However, symmetrically separated Lagrangian boundary conditions are specific Lagrangian boundary conditions.
\end{remark}

\begin{definition}[pitchfork bifurcation]\label{def:pitchfork}
The set
\[
\{ (\mu,x)\, | \, \nabla_x F(\mu,x) =0\},
\]
is a catastrophe set of a pitchfork bifurcation if and only if there exists a morphism from the scalar valued map $F$ to $f(\mu,x) = x^4+\mu x^2$ in the category of unfoldings. 
\end{definition}

Roughly speaking, a pitchfork bifurcation has a catastrophe set looking qualitatively like figure \ref{fig:HamPitchInt} or figure \ref{fig:bifur_pitchfork_time_reversal}. See \cite[Ch. 3]{lu1976singularity} for a definition of the category of unfoldings in the framework of catastrophe theory with respect to right and right-left equivalence. Notice that $f$ is not a versal unfolding of the map germ $x^4$ in the usual catastrophe theory framework (right- or right-left equivalence)
since there is no morphism from the universal unfolding of the cusp $x^4+\mu_2 x^2 +\mu_1 x$ to the unfolding $f$ from definition \ref{def:pitchfork}. This is the reason why the pitchfork bifurcation is not a generic phenomenon in gradient-zero-problems. However, in the presence of a $\Z / 2\Z$-symmetry, it is generic \cite[p.352]{poston1978catastrophe} \cite[p.486]{wassermann1988classification}. In Arnold's classification it occurs as the singularity $B_2$ \cite{ArnoldSym}. Indeed, in \ref{subsubsec:symexplperiodicpitchfork} we will show that a $\Z / 2\Z$-symmetry is present.

\begin{theorem}[generic occurrence of pitchfork bifurcations in integrable systems]\label{thm:genericpitchfork}
Consider a one-parameter family of $2n$-dimensional completely integrable Hamiltonian systems with time-$\tau$-flow maps $\phi_\mu$ and symmetrically separated boundary conditions for $\phi_\mu$ defined by Lagrangian manifolds $G_\mu$. Assume that for $\mu=0$ a compact, resonant Liouville torus $T$ (common level set of the integrals of motion) with period $\tau$ intersects $G_0$ in an isolated point $p$ such that $T_p T \cap T_p G$ is one-dimensional. Then either the catastrophe set of the boundary value problem shows a pitchfork bifurcation at $(\mu,z)=(0,p)$ or the problem is degenerate.
\end{theorem}


\textit{Proof.}
Scaling the Hamiltonian $H_\mu$, we can assume that $\tau = 1$.
Restricting to a neighbourhood of $T$ and $\mu=0$, there exist action angle coordinates
$(\theta_\mu,I_\mu)= ({\theta_\mu}_1,\ldots, {\theta_\mu}_n,I_\mu^1,\ldots,I_\mu^n)$ such that
$\theta_\mu(p)=0$, $I_\mu(p)=0$
and $I_\mu^1,\ldots,I_\mu^n$ are constants of motion. Viewing the systems in these coordinates, we do not need to keep track of the initial $\mu$-dependence and denote the coordinates as $(\theta,I)$.
Let us refer to the $\theta$-component of the symplectic map $\phi_\mu$ as $\phi^\theta_\mu$.\\

It suffices to consider a constant Lagrangian family $G_\mu=G$. Assume that the submanifold $G$ can (locally) be parametrised with $\theta$. Since $G$ is Lagrangian, there exists a scalar function $g$ with $g(0)=0$ such that $G=\{ (\theta,\nabla g(\theta)) \, | \, \theta \in (\R / 2\pi\Z)^{n} \}$. The intersection of the tangent space $T_pT$ and $T_pG$ is one-dimensional such that the Hessian matrix of $g$ at $p$ has exactly one vanishing eigenvalue. By the Splitting Lemma \cite[Thm.\ 5.3]{lu1976singularity}, locally around $\theta=0$ there exists a change in the $\theta$-coordinates fixing $0$ such that
\[
g(\theta)=h(\theta_1)+q(\theta_2,\ldots,\theta_n), \quad \text{with}\quad h(0)=h'(0)=h''(0)=0
\]
and $q(\theta_2,\ldots,\theta_n) = \sum_{j=2}^n \e_j \theta_j^2$ is a non-degenerate quadratic map with $\e_j \in \{-1,1\}$. Define $\mathrm{Sgns} := \mathrm{diag}(\e_1,\ldots,\e_n)$ for later.\\

\begin{figure}
\includegraphics[width=0.3\textwidth]{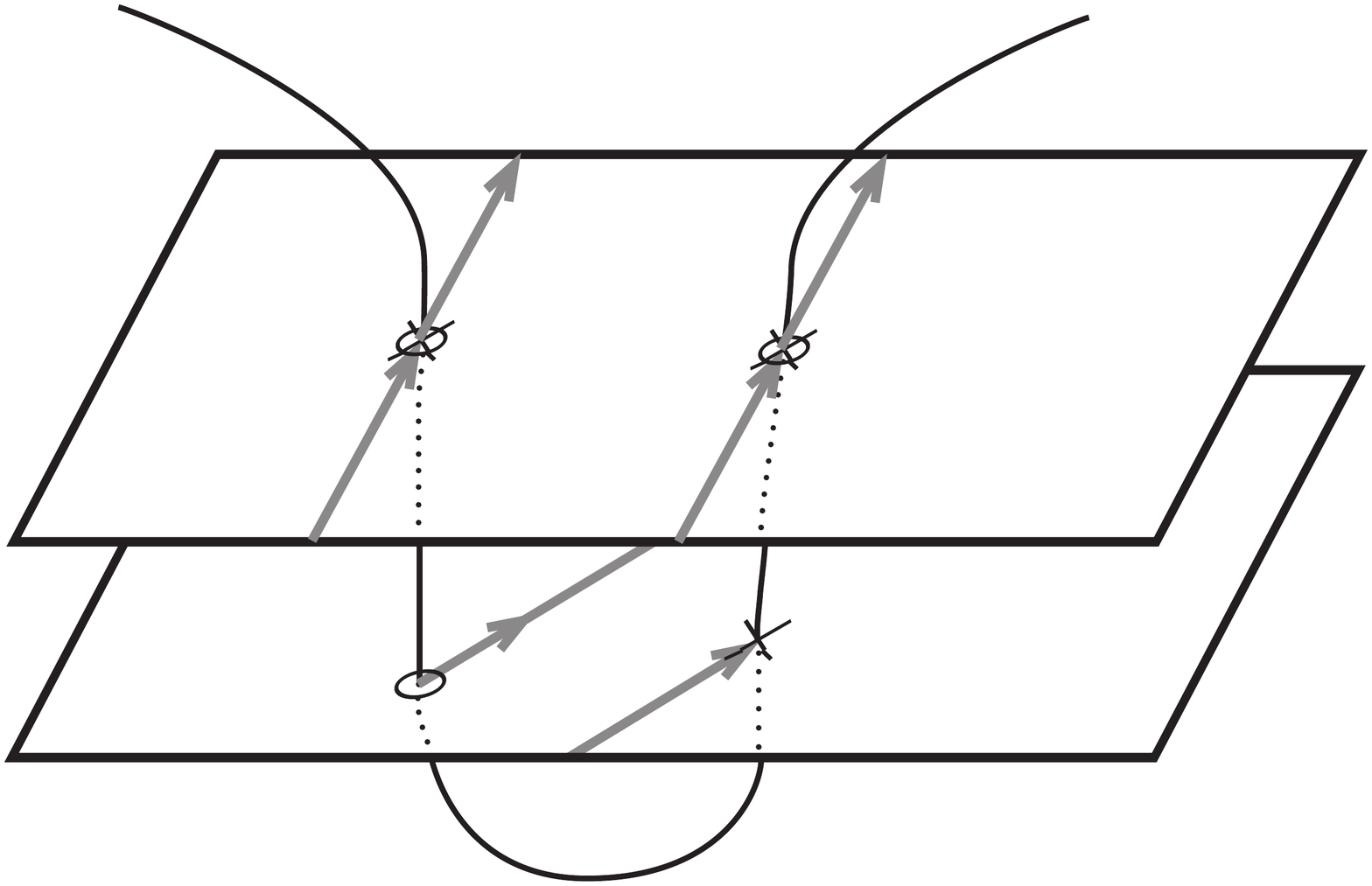}
\includegraphics[width=0.3\textwidth]{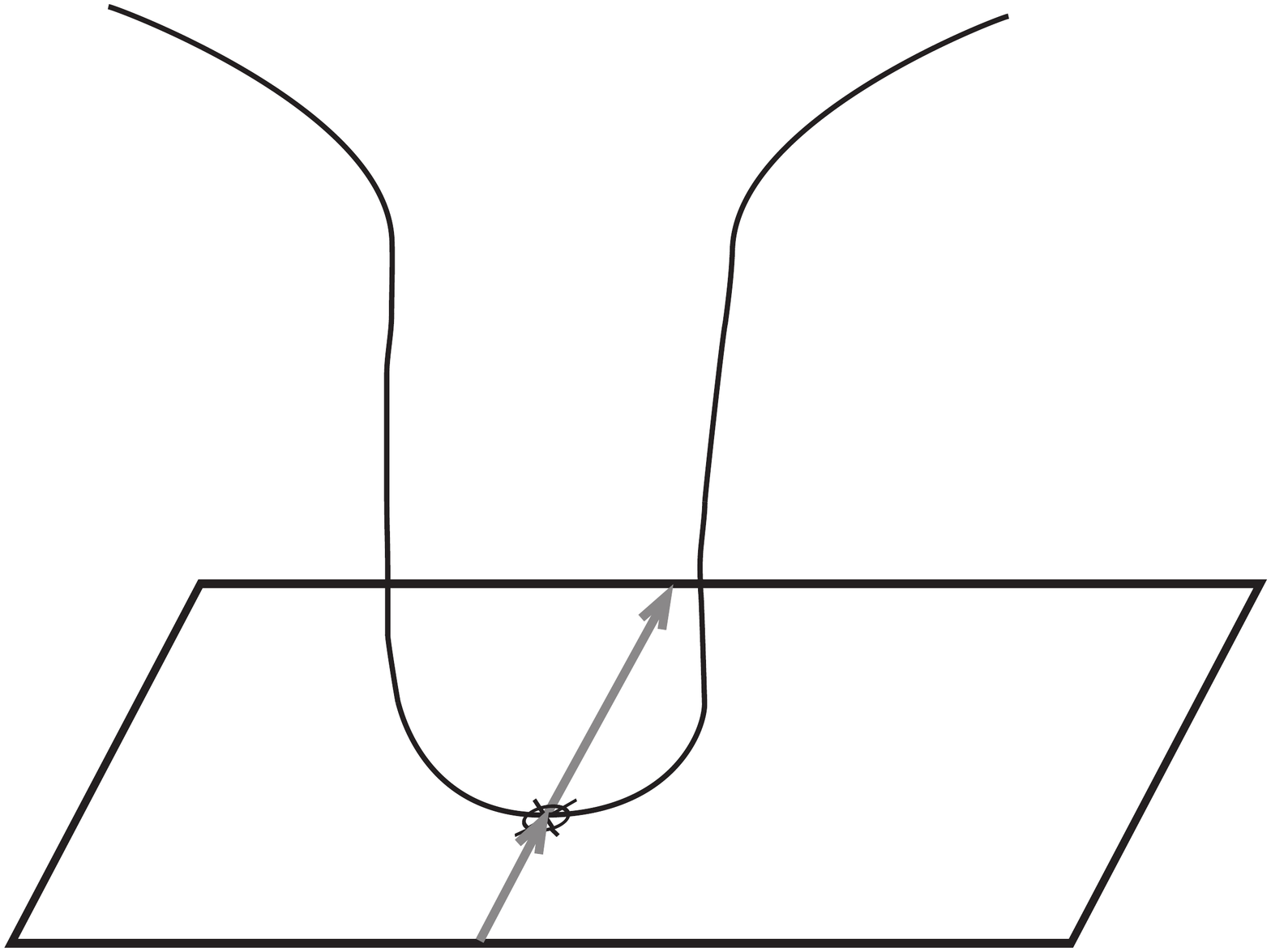}
\includegraphics[width=0.3\textwidth]{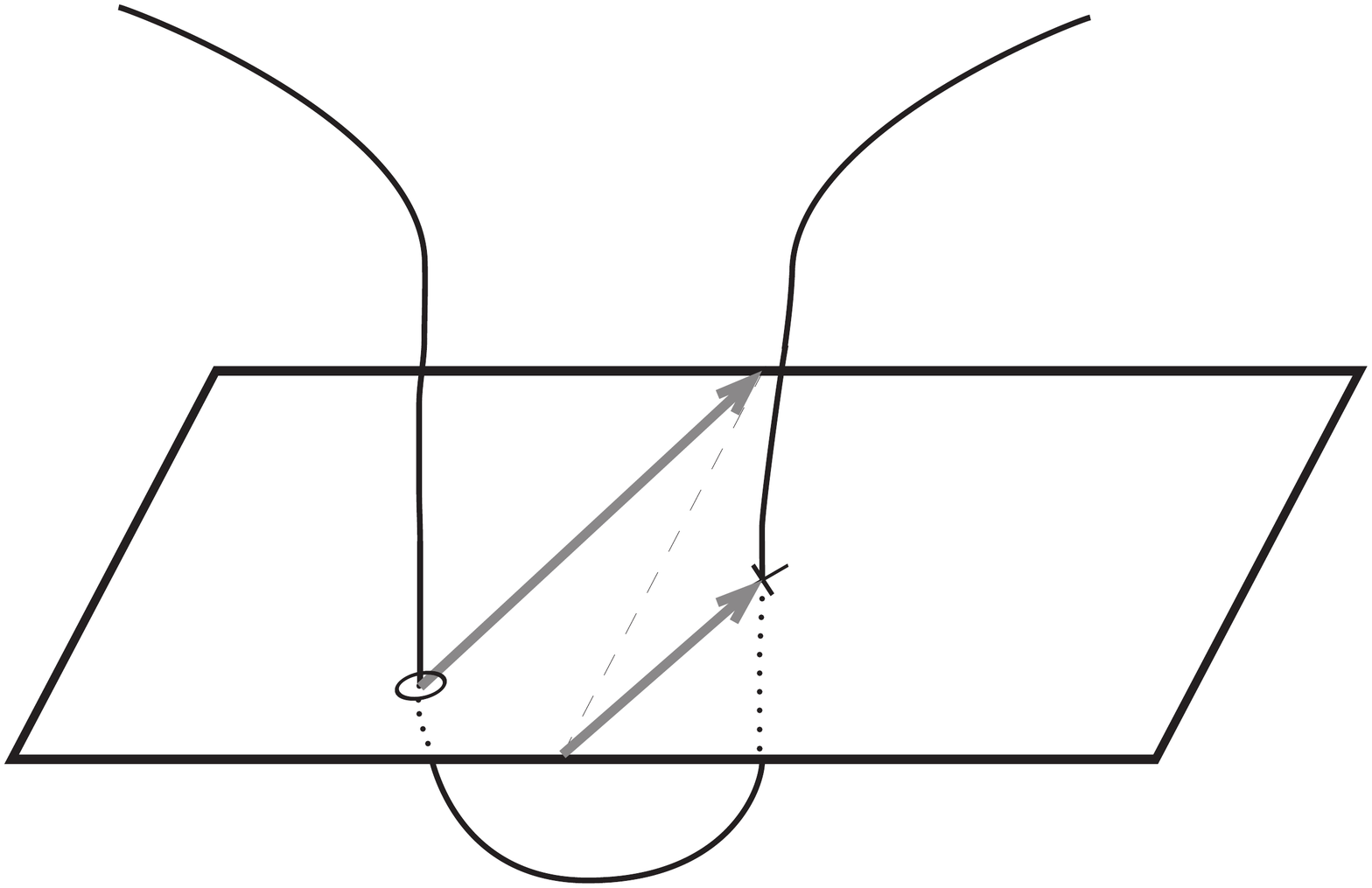}
\caption{Illustration of the mechanism in a four-dimensional phase space before, at and after passing through the bifurcation point. Opposite edges of the parallelograms are identified. They represent Liouville tori. The line intersecting the Liouville tori represents the manifold $G_\mu$. Gray arrows represent motions of the system starting and ending at an intersection point of a Liouville torus with $G_\mu$ within time $\tau=1$. }\label{fig:mechanismbifur}
\end{figure}

A motion starting at an intersection point $(\theta,I) = (\theta,\nabla g(\theta))$ of a Liouville torus with $G$ is a solution to the boundary value problem if and only if the endpoint $(\phi_\mu^\theta (\theta,I),I)$ of the motion also lies on $G$. Therefore, $(\theta,I)$ solves the boundary value problem if and only if $I=\nabla g(\theta)$ and $\theta$ fulfils
\begin{align*}
0 = b_\mu(\theta)
 &= \nabla g \left( \phi_\mu^\theta (\theta,\nabla g(\theta))\right) - \nabla g(\theta) .
\end{align*}
By the assumptions on the intersection of $T$ and $G$ 
\begin{itemize}\label{list:bullylist}
\item
the Liouville torus $T$ intersects with $G$ in $(\theta,I)=(0,0)$, i.e.\ $\nabla g(0) = 0$.
\item
The intersection of $T$ and $G$ is tangential such that $h(0)=h'(0)=h''(0)=0$ as obtained by the Splitting Lemma.
\item
All motions on the Liouville torus $T$ are periodic with period 1, i.e.\ $\phi_0^\theta(\theta,0) = \theta$ such that the Jacobian $\D_\theta \phi_0(\theta,0)$ is the identity matrix.
\end{itemize}
The situation is illustrated in the schematic picture in the middle of figure \ref{fig:mechanismbifur}.
Locally around $\theta=0$, the first component of $b_\mu$ is given as
\[
b^1_\mu(\theta) 
= h'\left( \phi^{\theta_1}_\mu (\theta,\nabla_\theta g(\theta)) \right)-h'(\theta_1). 
\]
Using the statements formulated in the bullet point list, a Taylor series expansion of $b^1_\mu(\theta)$ around $\mu=0$, $\theta=0$ has no constant term and the coefficients of $\mu$, $\theta_1,\ldots,\theta_n$ and $\theta_1^2$ vanish. The other coefficients are non-zero under non-degeneracy assumptions on the problem. Moreover, for the remaining components $b_\mu^{2,\ldots,n}(\theta)$ the Jacobian $D_{\theta_2,\ldots,\theta_n} b_0^{2,\ldots,n}(0) = \mathrm{Sgns} \cdot D_{I_2,\ldots,I_n} \phi_0^{\theta_2,\ldots,\theta_n}(0,0)$ has full rank (again by non-degeneracy). By the implicit function theorem, locally around $(\mu,\theta)=(0,0)$ there exist functions $\theta_j(\mu,\theta_1)$ with $\theta_j(0,0)=0$ ($j\ge 2$) such that 
\[ b_\mu^{2,\ldots,n}(\theta_1,\theta_2(\mu,\theta_1),\ldots,\theta_n(\mu,\theta_1))=0.\]

Moreover, using $\left.\frac{\p b^j}{\p \theta_1}\right|_{(\mu,\theta)=(0,0)}=0$ and $\left.\frac{\p b^j}{\p \theta_i}\right|_{(\mu,\theta)=(0,0)}\not=0$ for $i,j\ge2$, we conclude $\frac{\p \theta_j(0,\theta_1)}{\p \theta_1}(0)=0$. Therefore, in the Taylor series expansion of
\[ b_\mu^1(\theta_1,\theta_2(\mu,\theta_1),\ldots,\theta_n(\mu,\theta_1))\]
around $\mu=0$ and $\theta=0$ there is (still) no constant term and the coefficients of $\mu$, $\theta_1$ and $\theta_1^2$ vanish. Thus, there is a pitchfork bifurcation at $(\mu,\theta,I)=(0,0,0)$. 
\qed\\

We conclude that the pitchfork bifurcation occurs generically in symmetrically separated 1-parameter boundary value problems in completely integrable Hamiltonian systems. The mechanism of the bifurcation in the phase space is illustrated in figure \ref{fig:mechanismbifur}.\\


\begin{remark} The map $h$ in the proof of theorem \ref{thm:genericpitchfork} is a map germ $(\R,0)\to (\R,0)$ with
\[h(0)=h'(0)=h''(0)=0\]
and, under non-degeneracy assumptions, $h'''(0)\not=0$. This means that $h$ has a singularity of type $A_2$ (fold). A fold singularity is a generic phenomena in one-parameter families of smooth maps.
Allowing more parameters and arbitrarily high dimensional phase spaces, we can achieve any local Lagrangian intersection problem as an intersection of $G_\mu$ and $T$. Therefore, in a generic setting, $h$ can have any singularity that occurs for generic $k$ parameter families of smooth maps in the gradient-zero-problem. The type of singularity corresponds to the type of contact of the boundary condition and the Liouville torus at the bifurcation point. Thus, an analysis analogous to the proof of theorem \ref{thm:genericpitchfork} can be carried out for the other normal forms from Arnold's ADE-classification.
\end{remark}

%
%
%



\subsubsection{Symmetry based explanation of the periodic pitchfork}\label{subsubsec:symexplperiodicpitchfork}

Let us view the periodic pitchfork bifurcation as an effect of symmetry present in completely integrable systems. Indeed, we will show how complete integrability leads to a local $\Z / 2 \Z$ mirror symmetry in the generating function for the Hamiltonian diffeomorphism related to the system. Thus, the periodic pitchfork is connected to the bifurcation $B_2$ in Arnold's list for the $\Z / 2\Z$-symmetric gradient-zero problem \cite{ArnoldSym}.
Other references for the gradient-zero-problem with symmetry are \cite{wassermann1988classification} and \cite[Ch.14 \S 15]{poston1978catastrophe}. The following may be seen as an alternative proof of theorem \ref{thm:genericpitchfork}.\\

Consider a one-parameter family of completely integrable Hamiltonian systems of fixed dimension $2n$ with Hamiltonians $H_\mu$. Consider symmetrically separated Lagrangian boundary conditions. Suppressing the $\mu$-dependence of the action angle coordinates, the equations of motions read
\begin{align}\nonumber
\dot I &= 0\\ \nonumber
\dot \alpha &= \nabla_I H_\mu(I). \label{eq:aacords}
\end{align}
The time-1-map of the flow maps $(I,\alpha)$ to $(K,\beta)$ with
\begin{align}
K &= I\\ \nonumber
\beta &= \nabla_I H_\mu(I) + \alpha. \label{eq:aamotion}
\end{align}
Let us assume that $\nabla_I H_0(0)=0$. Under non-degeneracy conditions on the family $H_\mu$, locally near $I=0$ and $\mu=0$ there exists an inverse of the map $\nabla_I H_\mu$ fixing $0$. Let us denote its local inverse by $(\nabla_I H_\mu)^{-1}$ and its primitive by $\mathcal H_\mu$. The time-1-map can now be obtained by the generating function $S_\mu(\alpha,\beta) = \mathcal H_\mu(\beta-\alpha)$.\\
 
Assume that locally near $(I,\alpha)=(0,0)$, the boundary condition can be expressed as $I=\tilde g_\mu(\alpha)$, $K=\tilde g_\mu(\beta)$.
Since the boundary condition is Lagrangian, there exists a local map $g$ with $g(0)=0$ such that $\nabla g(\alpha)=\tilde g(\alpha)$. Assume that there exists a local change of coordinates in the $\alpha$ variables fixing 0 such that $g$ is of the form
\[
g(\alpha) = \alpha_1^3+q(\alpha_2,\ldots,\alpha_n).
\]
(Recall from the proof of theorem \ref{thm:genericpitchfork} that this holds generically at a 1-dimensional touch of the boundary condition with an invariant submanifold.) The change of variables in $\alpha$ can be extended to a symplectic change of variables $(I,\alpha)$ in the phase space fixing $(I,\alpha)=(0,0)$. 
The equations of motion keep their structure and we denote the new coordinates and maps by the same symbols as before. 
The boundary condition can be obtained from the generating function $B(\alpha,\beta )=\alpha_1^3-\beta_1^3+q(\alpha_2,\ldots,\alpha_n)-q(\beta_2,\ldots,\beta_n)$. Solutions to the boundary value problem correspond to critical points of
\[
F_\mu(\alpha,\beta)=B(\alpha,\beta )-\mathcal H_\mu(\beta-\alpha)=\alpha_1^3-\beta_1^3-\mathcal H_\mu(\beta-\alpha)+q(\alpha_2,\ldots,\alpha_n)-q(\beta_2,\ldots,\beta_n).
\]
Notice that $F_\mu$ is invariant under the transformation $(\alpha_1,\beta_1)\mapsto (-\beta_1,-\alpha_1)$. Now, by Arnold's classification \cite{ArnoldSym} bifurcations on the fix point set of the transformation correspond to $\Z/2 \Z$ bifurcations. Under the given assumptions, the bifurcation $B_2$ occurs at $(I,\alpha)=(0,0)$.\qed

\subsection{Effects of symmetry}\label{subsec:symeffects}

In the previous section \ref{subsec:completeIntegrable} it is shown that in boundary value problems extra structure (complete integrability) can lead to the occurrence of extra bifurcations which are not generic in the gradient-zero-problem.
Clearly, if in the Singularity-Lagrangian linkage theorem \ref{thm:linkageThm} the generating functions $f_\mu$ of symplectic maps and the generating functions $g_\mu$ of the boundary conditions obey the same symmetry relation then the bifurcations of the boundary value problem will be governed by the gradient-zero-problem for scalar valued maps obeying that symmetry. In the following we will analyse how symmetries of Hamiltonian systems, (or, more generally, of family of symplectic maps and their boundary conditions), translate to symmetries of generating functions and are, therefore, relevant to predict and explain bifurcation behaviour.

\subsubsection{Classical time reversal symmetry}\label{subsec:timereversalsym}
In view of the importance of classical Hamiltonian mechanics, we first show a numerical example illustrating the mechanism of a pitchfork bifurcation due to time reversal symmetry and then explain the general impact of time-reversibility on mechanical systems in their standard form. In section \ref{subsubsec:symmoregeneral} we proceed to a more general treatment of symmetries and reversal symmetries.

\paragraph{Numerical example}
Consider the time reversal Hamiltonian system $\R^2$ with the symplectic form $\d x \wedge \d y$ and Hamiltonian
\[
H(x,y)=\cos(y^2)+\mu x^2+x^3.
\]
We consider the Dirichlet boundary value problem $\phi^X(1,y)=1$ where $\phi$ is the time-$0.1$-map of the Hamiltonian flow.
In other words, a motion of the Hamiltonian system solves the boundary value problem if and only if it starts and ends after time $\tau=0.1$ on the line $\{1\}\times \R$ in the phase space. The described boundary condition is Lagrangian, symmetrically separated (definition \ref{def:symbd}) and $\Z / 2\Z$-symmetric w.r.t.\ the $x$-axis.\\

\begin{figure}
\begin{center}
\includegraphics[width=0.7\textwidth]{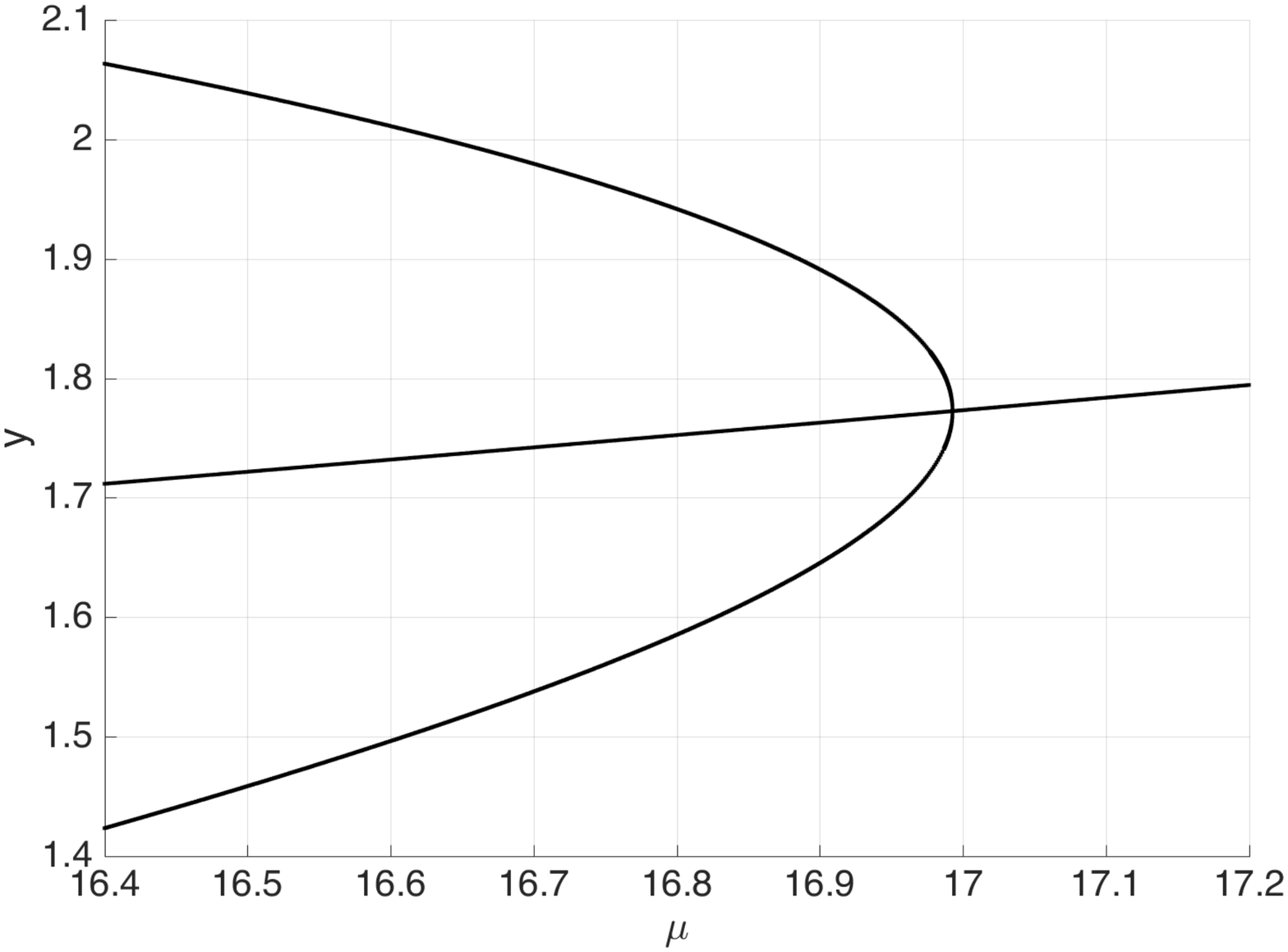}
\end{center}
\caption{Bifurcation diagram for a time-reversal Hamiltonian and boundary condition $x(0)=1=x(0.1)$. The Hamiltonian flow is obtained using the symplectic leapfrog method with step size $0.0005$.}\label{fig:bifur_pitchfork_time_reversal}
\end{figure}

\begin{figure}
\begin{center}
\includegraphics[width=0.49\textwidth]{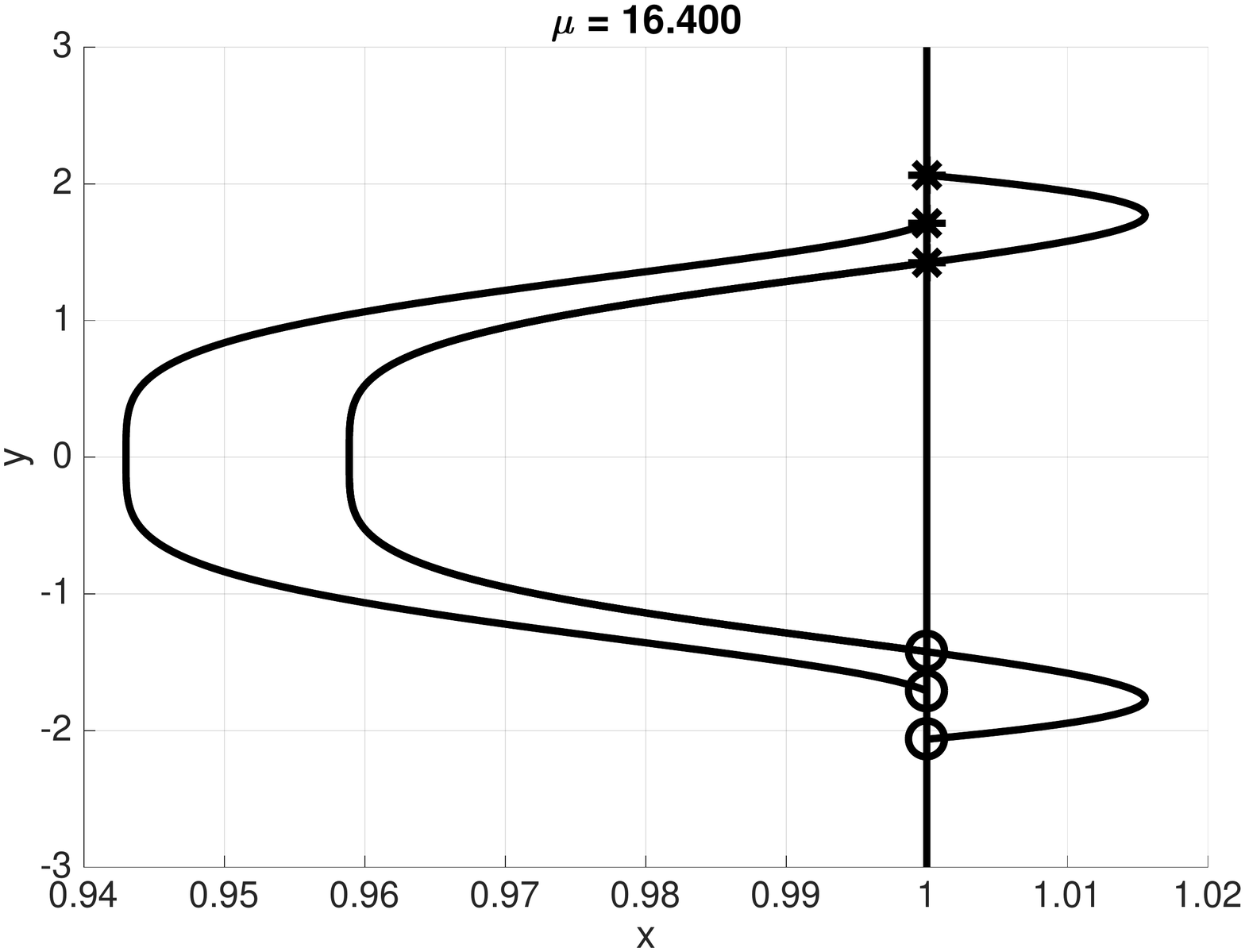}
\includegraphics[width=0.49\textwidth]{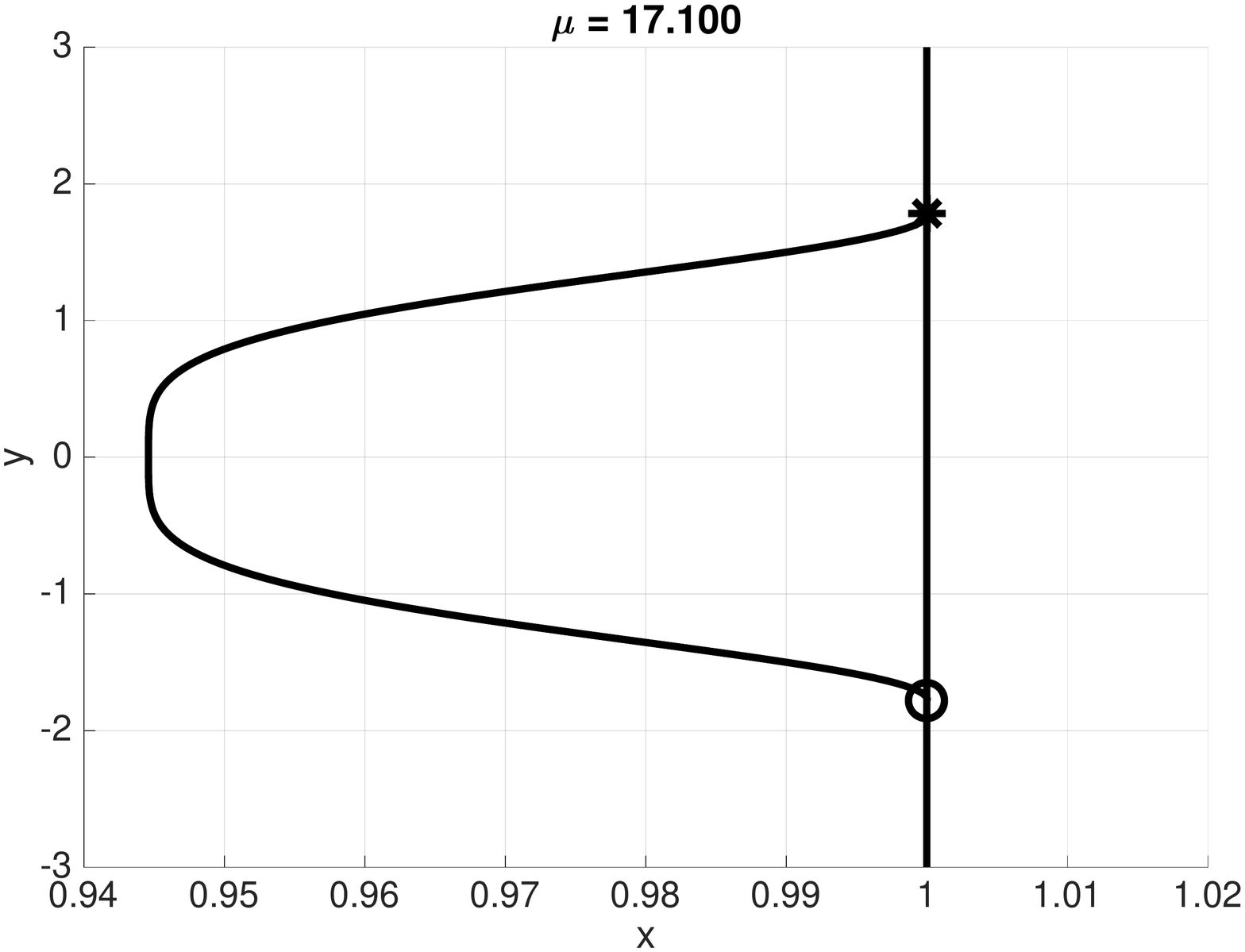}
\end{center}
\caption{Motions involved in the pitchfork bifurcation shown in figure \ref{fig:bifur_pitchfork_time_reversal}.} \label{fig:motion_pitchfork_time_reversal}
\end{figure}

Using the symplectic leapfrog method to integrate Hamilton's equations we obtain a bifurcation diagram showing a pitchfork bifurcation (figure \ref{fig:bifur_pitchfork_time_reversal}).
The motions involved in the bifurcation are shown in figure \ref{fig:motion_pitchfork_time_reversal} for two different values of $\mu$. We see that one of the orbits constitutes two solutions to the boundary value problem making use of the time reversal symmetry of the system. This orbit merges with another orbit that solves the boundary condition, exactly where both orbits become tangent to the boundary condition. At this point the relation $\frac{\p H_\mu}{\p y}(1,y)=0$ is necessarily fulfilled.\\

Indeed, the $\Z / 2\Z$ symmetry of the Hamiltonian induces a time reversal symmetry of the flow which the boundary condition respects such that the bifurcations of the system are governed by the gradient-zero-problem with symmetry. Let us describe the symmetry in a more general setting in terms of generating functions.

\paragraph{Mechanical Hamiltonians}\label{p:timereversalMechHam}
Let $(x,y)=(x^1,\ldots,x^n,y_1,\ldots,y_n)$ be symplectic coordinates and $H_\mu(x,y)$ a family of Hamiltonian functions such that $H_\mu(x,y)=H_\mu(x,-y)$. For example, mechanical Hamiltonians $H_\mu(x,y)=\frac 12\langle y,y\rangle - V(x)$, where $\langle . , . \rangle$ is the euclidean scalar product and $V$ a scalar valued map, fulfil this condition.
Denote the time-$\tau$-map of the Hamiltonian system $H_\mu$ by $\phi_\mu$ and define $X_\mu=x\circ \phi_\mu$ and $Y_\mu=y\circ \phi_\mu$. By the symmetry of $H_\mu$, the Hamiltonian system is time reversal symmetric such that
\begin{equation}\label{eq:time-rev-phi}
 \phi_\mu\big(X_\mu(a,b),-Y_\mu(a,b)\big)=(a,-b)
\end{equation}
for all points $(a,b)$ in the phase space. 
Assume $\det \left(\frac{\p Y_\mu}{\p x}\right)_{i,j} \not =0$ such that $y,Y_\mu$ constitutes a coordinate system. There exist generating functions $S_\mu(y,Y)$ such that
\[
\begin{pmatrix}x \\ -X_\mu\end{pmatrix}
=\nabla S_\mu(y,Y_\mu).
\]
Let $\xi(y,Y)=(-Y,-y)$. By \eqref{eq:time-rev-phi} the following holds:
\[
\nabla S_\mu = \begin{pmatrix}
x\\-X_\mu
\end{pmatrix}= \begin{pmatrix}
X_\mu \circ \xi \\ -x \circ \xi
\end{pmatrix} = \nabla (S_\mu \circ \xi).
\]
In this way, the symmetry of the problem is reflected in the corresponding generating functions. Now, if a boundary condition for $\phi_\mu$ can be represented with a generating function $B_\mu(y,Y)$ of the same type as $S_\mu$ such that $\nabla B_\mu = \nabla (B_\mu \circ \xi)$ then the gradient-zero-problem for $S_\mu-B_\mu$ will be governed by singularity theory for maps with the symmetry $\xi$.\\

If, for example, the boundary conditions are symmetrically separated with Lagrangian manifold $G_\mu$ in the phase space and if $y$ constitutes a coordinate system for $G_\mu$  (as in the numerical example), then $G_\mu = \{(y,\nabla b_\mu (y))\}_y$ for some scalar valued map $b_\mu$. We can choose $B_\mu(y,Y)=b_\mu(y)-b_\mu(Y)$. Then $\nabla B_\mu = \nabla (B_\mu\circ \xi)$. In the numerical example we have $G_\mu=\{1\}\times \R$, $b_\mu(y)=y$ and $B_\mu(y,Y)=y-Y$. \\


We conclude that extra symmetry in Lagrangian boundary value problems can allow for the generic occurrence of bifurcations, which are non-generic in the class of Lagrangian boundary value problems, if the symplectic map and the boundary conditions obey the same symmetry relation. The problem reduces to the gradient-zero-problem with symmetry. In the following section \ref{subsubsec:symmoregeneral} we provide a more general treatment of symmetries and reversing symmetries on bifurcation problems.

\subsubsection{General treatment of symmetries and reversing symmetries}\label{subsubsec:symmoregeneral}

Recall that it is sufficient to consider locally defined, symplectic maps $\phi$ which map 0 to 0 in the space $V = U \oplus U^\ast$ as we can locally trivialize symplectic manifolds $M$ around each $z\in M$ and $\phi(z) \in M$ using centred Darboux coordinates. In the following we will use the setting of section \ref{subsub:coordfree} to describe the effects of ordinary and reversing symmetries on generating functions and, thus, on the local bifurcation behaviour of Lagrangian boundary value problems. In order to avoid cumbersome notation and repeating remarks that domains of definition might shrink, we neglect to incorporate in our notation that maps are defined only on neighbourhoods of 0 of the occurring spaces.\\

Consider a symplectic map $\phi \colon V=U \oplus U^\ast \to \underline V=\underline U \oplus \underline U^\ast$ mapping 0 to 0 and a diffeomorphism $\Psi\colon V \to V$. Let $\Psi^1$ denote the $U$-component of the map $\Psi$ and $\Psi^2$ the $U^\ast$-component, i.e.\ $\Psi(v)=(\Psi^1(v),\Psi^2(v)) \in U \oplus U^\ast$. Moreover, assume there exists a generating function $S \colon U \oplus \underline U \to \R$ defined around 0 such that for all $(u,\underline u)\in U \oplus \underline U$ the following Lagrangian submanifolds of $(V \oplus \underline V, \omega \oplus (-\omega))$ coincide:
\[
\Gamma := \{ (v,\phi(v)) | v \in V \} = \{\big(u, \D_uS(u,\underline u), \underline u, -\D_{\underline u}S(u,\underline u)\big)\,|\, (u,\underline u) \in U \oplus \underline U\}.
\]

\paragraph{Ordinary symmetry}

\begin{prop}\label{prop:ordinsym} The relation
\begin{equation}\label{eq:ordinarysymmetry}
\Psi^{-1} \circ \phi \circ \Psi = \phi
\end{equation}
is equivalent to $\forall (u,\underline u) \in U \oplus \underline U :$
\begin{align}\nonumber
\D_u S\Big( \Psi^1(u,\D_uS(u,\underline u)),\Psi^1(\underline u,-\D_{\underline u}S(u,\underline u))\Big) &= \Psi^2\Big(u, D_u S(u,\underline u)\Big)\\ \label{eq:ordinsymRel}
-\D_{\underline u} S\Big( \Psi^1(u,\D_uS(u,\underline u)),\Psi^1(\underline u,-\D_{\underline u}S(u,\underline u))\Big) &= \Psi^2\Big(\underline u, - D_{\underline u} S(u,\underline u)\Big)
\end{align}
\end{prop}

\begin{remark} If $\Psi$ is additionally symplectic then $\Psi$ is referred to as an \textit{ordinary symmetry} for $\phi$.
If, for instance, the map $\phi$ is given as a Hamiltonian diffeomorphism and $\Psi$ is a symplectic map on the phase space leaving the Hamiltonian invariant then $\Psi$ as well as $\Psi^{-1}$ are ordinary symmetries for $\phi$.
\end{remark}

\textit{Proof of proposition \ref{prop:ordinsym}.} The relation $\phi \circ \Psi=\Psi \circ \phi$ is equivalent to $(\Psi \times \Psi)(\Gamma)=\Gamma$, where $(\Psi \times \Psi)\colon V \oplus \underline V \to V \oplus \underline V$, $(\Psi \times \Psi)(v,\underline v)=(\Psi(v),\Psi(\underline v))$ denotes the diagonal action. The relation $(\Psi \times \Psi)(\Gamma)=\Gamma$ is equivalent to \eqref{eq:ordinsymRel}.\qed\\

The following proposition analyses the effects of symmetries on the phase space $U \oplus U^\ast$ which restrict to the Lagrangian submanifold $U$. Examples are spatial symmetries of mechanical systems.

\begin{prop}\label{prop:hdiagsym} Let $h \colon U \to U$ be a diffeomorphism. The symplectic map $\Psi = (\Psi^1,\Psi^2) \colon U \oplus U^\ast \to U \oplus U^\ast$ defined by
\begin{equation}\label{eq:defPsih}
\Psi^1(u,u^\ast)=h(u) \quad \Psi^2(u,u^\ast)=u^\ast \circ \D h^{-1}(u)
\end{equation}
is an ordinary symmetry for $\phi$ if and only if S is invariant under $h$, i.e.\
\[
S \circ (h \times h) = S,
\]
holds true (up to a constant), where $(h \times h)(u,\underline u)=(h(u),h(\underline u))\in U \oplus \underline U$ denotes the diagonal action.
\end{prop}

\begin{remark}
The symplectic map $\Psi$ defined in \eqref{eq:defPsih} corresponds to the cotangent lifted action of $h$. The reader might be familiar with its representation in Darboux coordinates $p,q$ of the cotangent bundle: $q \mapsto h(q), p \mapsto \D h(q)^{-t}p$. Moreover, we have used the cotangent lifted action in the proof of proposition \ref{prop:rightequivalence}.
\end{remark}

\textit{Proof of proposition \ref{prop:hdiagsym}.} The relation $S \circ (h \times h) = S$ (up to a constant) is equivalent to
\[
\D (S(h(u),h(\underline u)))=\D S(u,\underline u) \quad \forall (u,\underline u) \in U \oplus \underline U,
\]
which is equaivalent to $\forall (u,\underline u) \in U \oplus \underline U:$
\begin{align*}
\D_u (S(h(u),h(\underline u)))&=\D_u S(u,\underline u)\circ \D h^{-1}(u)\\
-\D_{\underline u} (S(h(u),h(\underline u))) &= - \D_{\underline u} S(u,\underline u)\circ \D h^{-1}(\underline u),
\end{align*}
which is equivalent to \eqref{eq:ordinsymRel} with $\Psi$ defined as in \eqref{eq:defPsih}.
The claim follows by proposition \ref{prop:ordinsym}.\qed\\

Analogous to lemma \ref{lem:graphoverzero}, after a symplectic change of coordinates in $V \oplus \underline V$, we can assume that boundary conditions as well as the graph $\Gamma$ of a symplectic map $\phi$ are both graphical over $U \oplus \underline U$.
We state the following proposition which gives a criterion when $\Psi$ is a symmetry for separated, Lagrangian boundary conditions.

\begin{prop}\label{prop:Lagrangiansym}
Let $B(u,\underline u)=b(u)-\underline b(\underline u)$ for local maps $b\colon U \to \R$, $\underline b\colon \underline U \to \R$. Consider
\[
G = \{(u, \D b(u))\,|\, u \in U\},  \quad
\underline G = \{(\underline u, \D \underline b(\underline u))\, |\, \underline u \in \underline U\},
\quad \Lambda := G \times \underline G \subset (U \oplus U^\ast) \oplus (\underline U \oplus \underline U^\ast)
\]
The relations \eqref{eq:ordinsymRel} stated with $B$ instead of $S$ are fulfilled 
if and only if
\begin{equation}\label{eq:Ginvariance}
\Psi(G)=G \quad \Psi(\underline G)=\underline G.
\end{equation}
\end{prop}

\textit{Proof.} Notice that $\Lambda$ is generated by $B$, i.e.\
\[
\Lambda = G\times G=\{ (u, \D_u B(u,\underline u),\underline u,-\D_{\underline u}B(u, \underline u)\,|\, (u,\underline u)\in U \oplus \underline U\}.
\]
The relation \eqref{eq:Ginvariance} holds true if and only if the submanifold $\Lambda$ is invariant under the diagonal action $\Psi\times\Psi$ if and only if \eqref{eq:ordinsymRel} stated with $B$ instead of $S$ holds true. \qed

\begin{corollary} If $\Psi$ as in \eqref{eq:defPsih} is an ordinary symmetry for a family of symplectic maps and a family of separated Lagrangian boundary conditions invariant under $\Psi$ is given, then the Lagrangian boundary value problem behaves like the gradient-zero problem $\tilde S(u,\underline u)=0$ with symmetry $\tilde S \circ (h\times h)=\tilde S$. Using the notation of the propositions, $\tilde S = S -B$.
\end{corollary}

\noindent\textbf{Examples. }

\begin{itemize}
\item
Spatial symmetries of mechanical systems on $U$ with phase space $U \oplus U^\ast$ can be phrased as maps $h \colon U \to U$. The symplectic map $\Psi$ defined as in \eqref{eq:defPsih} is an ordinary symmetry for the Hamiltonian diffeomorphism $\phi$. The generating function $S$ of $\phi$ fulfils $S \circ (h \times h)=S$.

\begin{itemize}
\item
Let $h(u)=-u$, $\Psi(u,u^\ast)=(-u,-u^\ast)$. It follows that $S$ fulfils the $\Z / 2 \Z$ symmetry relation $S(-u,-\underline u)=S(u,\underline u)$.

\item
For $\lambda \in \R\setminus\{0\}$ let $h_\lambda(u)=\lambda u$, $\Psi_\lambda(u,u^\ast)=(\lambda u, \lambda^{-1} u^\ast)$. Then $S(\lambda u, \lambda \underline u)=S(u,\underline u)$, i.e.\ $S$ is homogeneous (of degree 1).
\end{itemize}

\item
Consider the Dirichlet problem $u=0$, $\underline u=0$ for a symplectic map with generating function $S\colon U^\ast \oplus \underline U^\ast \to \R$ invariant under the cotangent lifted action of a linear, spatial symmetry $A \colon U \to U$. The problem behaves like the gradient-zero problem for $S$ with diagonal symmetry $A^t$, i.e.\ $S \circ (A^t \times A^t) = S$. Notice that the structure of the boundary condition forces us to use $(u^\ast,\underline u^\ast)$ as coordinates for the generating function causing the transposition of $A$ (compare to proposition \ref{prop:hdiagsym}).

\item
In case of $k$ independent, Poisson commuting integrals of motions, there exist symplectic coordinates such that the Hamiltonian $H \colon U\oplus U^\ast \to \R$ does not depend on $u^1,\ldots,u^k$. This means for each $\lambda \in \R^k$ the translation
\[
h(u^1,\ldots,u^k,u^{k+1},\ldots,u^n)=h(u^1+\lambda_1,\ldots,u^k+\lambda_k,u^{k+1},\ldots,u^n)
\]
gives rise to an ordinary symmetry for the Hamiltonian diffeomorphism. Its generating function $S\colon U \oplus \underline U \to \R$ depends on the following $2n-k$ variables: $u_1-\underline u_1,\ldots,u_k-\underline u_k,u_{k+1},\underline u_{k+1},\ldots,u_n,\underline u_{n}$. Notice that periodic boundary conditions
\[ \underline u - u =0, \quad
\underline u^\ast - u^\ast =0
\]
share this symmetry property. This recovers the fact that solutions to periodic boundary value problems for Hamiltonian diffeomorphisms with $k$ integrals of motion are $k$-dimensional manifolds. In planar Hamiltonian systems, for instance, periodic boundary value problems are solved by periodic orbits.

\end{itemize}

\paragraph{Reversal symmetry}

\begin{prop}\label{prop:reversalsym} The relation
\begin{equation}\label{eq:antisym} \Psi^{-1} \circ \phi \circ \Psi = \phi^{-1} \end{equation}
is equivalent to $\forall (u,\underline u)\in U \oplus \underline U:$
\begin{align}\nonumber
\D_u S\Big(\Psi^1(\underline u , - \D_{\underline u} S(u,\underline u)),\Psi^1( u ,  \D_{u} S(u,\underline u))\Big)
&= \Psi^2\Big(\underline u , - \D_{\underline u} S(u,\underline u)\Big)\\ \label{eq:antsymrel}
-\D_{\underline u} S\Big(\Psi^1(\underline u , - \D_{\underline u} S(u,\underline u)),\Psi^1( u ,  \D_{u} S(u,\underline u))\Big)
&= \Psi^2\Big(u ,  \D_{ u} S(u,\underline u)\Big).
\end{align}

\end{prop}

\begin{remark} If $\Psi$ is additionally anti-symplectic, i.e.\ $\Psi^\ast \omega= -\omega$, then $\Psi$ is referred to as a \textit{reversal symmetry} for $\phi$. For example, this situation arises if $\phi$ is given as a Hamiltonian diffeomorphism and $\Psi$ is a anti-symplectic map on the phase space leaving the Hamiltonian invariant. Inverting \eqref{eq:antisym}, it follows that $\Psi$ is a reversal symmetry for $\phi$ if and only if $\Psi^{-1}$ is a reversal symmetry.
\end{remark}

\textit{Proof of proposition \ref{prop:reversalsym}.} The relation $\phi \circ \Psi=\Psi \circ \phi^{-1}$ is equivalent to $(\Psi \underline \times \Psi)(\Gamma)=\Gamma$, where $(\Psi \underline \times \Psi)(v,\underline v)=(\Psi(\underline v), \Psi(v))$ denotes the reversal action. The invariance of $\Gamma$ is equivalent to \eqref{eq:antsymrel}.\qed\\

The following proposition analyses the effects of reversal symmetries on the phase space $V=U \oplus U^\ast$ which restrict to the Lagrangian submanifold $U$. Examples are spatial reversal symmetries.

\begin{prop}\label{prop:hdiagantisym} Let $h \colon U \to U$ be a diffeomorphism. The anti-symplectic map $\Psi = (\Psi^1,\Psi^2) \colon U \oplus U^\ast \to U \oplus U^\ast$ defined by
\begin{equation}\label{eq:defPsia}
\Psi^1(u,u^\ast)=h(u) \quad \Psi^2(u,u^\ast)=-u^\ast \circ \D h^{-1}(u).
\end{equation}
is a reversal symmetry for $\phi$ if and only if $S$ is invariant under $h$, i.e.\
\[
S \circ (h \underline \times h) = S,
\]
holds true (up to a constant), where $(h \underline \times h)(u,\underline u)=(h(\underline u),h(u))\in U \oplus \underline U$ denotes the reversal action.
\end{prop}

\textit{Proof.} The relation $S \circ (h \underline \times h) = S$ (up to a constant) is equivalent to
\[
\D (S(h(\underline u),h(u)))=\D S(u,\underline u) \quad \forall (u,\underline u) \in U \oplus \underline U,
\]
which is equaivalent to $\forall (u,\underline u) \in U \oplus \underline U:$
\begin{align*}
\D_{\underline u} (S(h(\underline u),h(u)))\circ \D h(u) &= \D_{u} S(u,\underline u)\\
\D_u (S(h(\underline u),h(u)))\circ \D h(\underline u)&=\D_{\underline u} S(u,\underline u)
\end{align*}
which is equivalent to \eqref{eq:antsymrel} with $\Psi$ defined as in \eqref{eq:defPsia}.
The claim follows by proposition \ref{prop:reversalsym}.\qed\\

The next proposition gives a criterion when an anti-symplectic map $\Psi$ is a reversal symmetry for separated, Lagrangian boundary conditions. It statement holds, however, for diffeomorphisms $\Psi \colon V \to V$.

\begin{prop}\label{prop:Lagrangianantisym}
Let $B(u,\underline u)=b(u)-\underline b(\underline u)$ for local maps $b\colon U \to \R$, $\underline b\colon \underline U \to \R$. Consider
\[
G = \{(u, \D b(u))\,|\, u \in U\},  \quad
\underline G = \{(\underline u, \D \underline b(\underline u))\, |\, \underline u \in \underline U\},
\quad \Lambda := G \times \underline G \subset (U \oplus U^\ast) \oplus (\underline U \oplus \underline U^\ast)
\]
The relations \eqref{eq:antsymrel} stated with $B$ instead of $S$ are fulfilled 
if and only if
\begin{equation}\label{eq:swapinv}
\Psi(G)=\underline G \quad \Psi(\underline G)= G.
\end{equation}
\end{prop}

\begin{remark} In contrast to proposition \ref{prop:Lagrangiansym}, the manifolds $G$ and $\underline G$ get swapped by the diffeomorphism $\Psi$.
\end{remark}

\textit{Proof of proposition \ref{prop:Lagrangianantisym}.}
Notice that $\Lambda$ is generated by $B$, i.e.\
\[
\Lambda = G\times G=\{ (u, \D_u B(u,\underline u),\underline u,-\D_{\underline u}B(u, \underline u)\,|\, (u,\underline u)\in U \oplus \underline U\}.
\]
The relation \eqref{eq:swapinv} holds true if and only if the submanifold $\Lambda$ is invariant under the reversal action $\Psi\underline\times\Psi$ if and only if \eqref{eq:antsymrel} stated with $B$ instead of $S$ holds true. \qed

\begin{corollary} If $\Psi$ as in \eqref{eq:defPsia} is a reversal symmetry for a family of symplectic maps and a family of separated Lagrangian boundary conditions invariant under $\Psi$ in the sense of \eqref{eq:swapinv} is given, then the local Lagrangian boundary value problem behaves like the gradient-zero problem $\tilde S(u,\underline u)=0$ with symmetry $\tilde S \circ (h\underline\times h)=\tilde S$. Using the notation of the propositions, $\tilde S = S -B$.
\end{corollary}

\textbf{Examples.} 
\begin{itemize}
\item
Let $h(u)=u$ such that $\Psi(u,u^\ast)=(u,-u^\ast)$. Then $S(u,\underline u)=S(\underline u,u)$. 

\item
Let $h(u)=-u$ such that $\Psi(u,u^\ast)=(-u,u^\ast)$. Then $S(u,\underline u)=S(-\underline u,-u)$. This recovers the results of section \ref{subsec:timereversalsym}.

\item
For $\lambda \in \R\setminus\{0\}$ let $h_\lambda(u)=\lambda u$, so $\Psi_\lambda(u,u^\ast)=(\lambda u, -\lambda^{-1} u^\ast)$. Then $S(\lambda \underline u, \lambda u)=S(u,\underline u)$. In particular $S(\underline u,u)=S(u,\underline u)$. The map $S$ is homogeneous and symmetric.

\end{itemize}

Let us relate these findings back to the periodic pitchfork described in \ref{subsubsec:periodicpitchforkexpl} and in \ref{subsubsec:symexplperiodicpitchfork}. Recall from the examples of ordinary symmetries that integrals of motions lead to translation symmetries. In completely integrable systems in action angle coordinates, the generating function $S \colon U \oplus \underline U \to \R$ of the Hamiltonian flow only depends on the difference $u-\underline u$ of the angles. In particular, $S$ is invariant under the transformation $(u_1,\underline u_1) \mapsto (- \underline u_1, -u_1)$.\\

Recall from \ref{subsubsec:symexplperiodicpitchfork} that under the assumption that a periodic orbit touches the boundary condition and gives rise to a solution of the boundary value problem, there exists a symplectic change of coordinates on $U \oplus U^\ast$ such that the symmetrically separated boundary condition is invariant under $(u_1,\underline u_1) \mapsto (- \underline u_1, -u_1)$ as well, locally around the touch point.
The bifurcation behaviour corresponds to the gradient-zero-problem with a $\Z / 2 \Z$ symmetry.
In a one-parameter family of systems with a $\Z / 2 \Z$ symmetry a pitchfork bifurcation is a generic phenomenon.\\

While the $\Z / 2 \Z$ symmetry in completely integrable systems is quite hidden and becomes visible using action angle coordinates, the $\Z / 2 \Z$ symmetry in time reversal systems is more apparent and again picked up by symmetrically separated Dirichlet boundary conditions since these share the phase space symmetry as elaborated in \ref{p:timereversalMechHam}. Again, a pitchfork bifurcation becomes a generic phenomena in one-parameter families of problems.

\section{Conclusion}

We have shown how Lagrangian boundary value problems for symplectic maps can be viewed as local intersection problems of Lagrangian submanifolds. This, in turn, corresponds to finding critical points of a smooth function (given as the difference of generating functions) and gives a {\em finite dimensional} approach to bifurcation theory for Hamiltonian boundary value problems.
This allows us to 
\begin{itemize}
\item
link generic problems to Thom's seven elementary catastrophes and to Arnold's ADE-classification, a framework known as catastrophe theory,
\item
explain (without genericity assumptions) how certain, typical boundary conditions can prohibit bifurcations,
\item
analyse how extra structure for symplectic maps can lead to extra bifurcations.
\end{itemize}

An instance of our analysis how boundary conditions can prohibit bifurcations is the following result: in Dirichlet problems for $2n$-dimensional systems only those bifurcations can occur which can be obtained in critical points of smooth function problems in at most $n$ variables. This means, for example, that in
Dirichlet problems for generic planar symplectic maps only $A$-series type bifurcations can occur.
In contrast, using periodic boundary conditions $D$-series singularities are also possible.
Moreover, extra structure for symplectic maps can lead to extra bifurcations. In symmetrically separated boundary value problems for flows of completely integrable Hamiltonian systems, we describe a pitchfork-type bifurcation, which we call the periodic pitchfork bifurcation. In this novel bifurcation, two complete periodic solutions bifurcate from a path of partial periodic orbits.\\



Lagrangian boundary value problems for symplectic maps correspond to gradient- zero-problems with symmetry if the symplectic map and the boundary values are governed by the same symmetry relation. Here, as well as generic bifurcations, bifurcations can occur which make use of the symmetry of the system. Propositions \ref{prop:Lagrangiansym} and \ref{prop:Lagrangianantisym} describe when separated Lagrangian boundary conditions fulfil a given symmetry. For instance, in one parameter families of time reversal symmetric Hamiltonian systems a time reversal pitchfork can occur generically.\\

For structurally simple symmetries, which split into separated actions on two Lagrangian subspaces (e.g. spatial symmetries), Propositions \ref{prop:hdiagsym} and \ref{prop:hdiagantisym} reveal which symmetry is induced in the corresponding gradient-zero-problem. In contrast, for arbitrary group actions we do not have a characterization of which symmetries induce a correspondence between boundary value problems and critical points of symmetric functions.\\

The framework presented here raises the possibility of discovering exotic bifurcations in examples from physics, and also of discovering further new phenomena induced in specific classes of equations.

\ack
This research was supported by the Marsden Fund of the Royal Society Te Ap\={a}rangi.

\section*{References}
\bibliography{resources}{}
\bibliographystyle{iopart-num} 


\end{document}